\documentclass[final]{siamltex}

\usepackage[latin1]{inputenc}
\usepackage[T1]{fontenc}
\usepackage[francais,english]{babel}
\usepackage{amsmath}
\usepackage{amssymb}
\usepackage{amssymb,amsmath,amsfonts,mathrsfs,graphicx}
\usepackage[colorlinks=true]{hyperref}

\newcommand{\Dh}{{\widehat{D}}}

\newcommand{\DDh}{{\widehat{d}}}
\newcommand{\yh}{{\widehat{y}}}
\newcommand{\Yh}{{\widehat{Y}}}
\newcommand{\Xh}{{\widehat{X}}}
\newcommand{\Lh}{{\widehat{\Lambda}}}
\newcommand{\pwh}{{\widehat{\bold{p_{\omega}}}}}
\newcommand{\Pwh}{{\widehat{\bold{P_{\omega}}}}}
\newcommand{\pvh}{{\widehat{\bold{p_v}}}}
\newcommand{\Pvh}{{\widehat{\bold{P_v}}}}
\newcommand{\pv}{\bold{p_v}}
\newcommand{\pw}{\bold{p_{\omega}}}
\newcommand{\ty}{\widetilde{y}}

\newcommand{\tX}{{\widetilde{X}}}
\newcommand{\tL}{{\widetilde{\Lambda}}}
\newcommand{\tD}{\widetilde{D}}
\newcommand{\tpv}{\widetilde{\bold{p_v}}}
\newcommand{\tpw}{\widetilde{\bold{p_{\omega}}}}
\newcommand{\SSS}{{\mathbb  S}}

\newcommand{\dv}[2]{{\frac{\partial #1}{\partial #2}}}

\newcommand{\dotex}{{\tfrac{d}{dt}}}

\newcommand{\R}{\mathbb{R}}

\newcommand{\NN}{{\mathbb N}}

\newcommand{\bu}{{\bold{u}}}
\newcommand{\bv}{{\bold{v}}}
\newcommand{\bw}{{\bold{w}}}
\newcommand{\wm}{\boldsymbol{\omega_m}}
\newcommand{\vm}{\boldsymbol{v_m}}
\newcommand{\om}{\boldsymbol{\omega}}
\newcommand{\vit}{\boldsymbol{v}}
\newcommand{\etab}{\boldsymbol{\eta}}

\title{Inertial-sensor bias estimation from  brightness/depth images and based  on $SO(3)$-invariant  integro/partial-differential equations on the unit sphere}
\author{
Nad\`ege Zarrouati-Vissi\`ere \thanks{DGA, 7-9 rue des Mathurins, 92220 Bagneux, FRANCE. Centre Automatique et Syst\`emes, Mines-ParisTech,
60, boulevard Saint-Michel
75272 Paris Cedex, FRANCE. email: nadege.zarrouati@mines-paristech.fr
}
\and
Karine Beauchard\thanks{
Centre de Math\'{e}matiques Laurent Schwartz, \'{E}cole Polytechnique, 91128 Palaiseau Cedex, FRANCE.
email: Karine.Beauchard@math.polytechnique.fr}
\and
Pierre Rouchon \thanks{Centre Automatique et Syst\`emes, Mines-ParisTech,
60, boulevard Saint-Michel
75272 Paris Cedex, France. email:pierre.rouchon@mines-paristech.fr
}
}

\date{September 2013}

\begin{document}

\maketitle

\begin{abstract}
Constant biases associated to measured linear and angular velocities of a moving object can be estimated from  measurements of a static  scene  by  embedded brightness  and depth sensors. We propose here a  Lyapunov-based observer taking  advantage of the $SO(3)$-invariance of the  partial differential equations satisfied  by the measured  brightness and  depth fields. The resulting asymptotic observer is governed by a non-linear integro/partial differential system where the two  independent scalar variables indexing the pixels  live on $\SSS^2$.  The observer design  and analysis are  strongly simplified by  coordinate-free  differential calculus on $\SSS^2$ equipped with its natural   Riemannian structure.  The observer convergence is investigated under  $C^1$ regularity assumptions  on the object motion and its scene. It relies on Ascoli-Arzela theorem and pre-compactness of the observer trajectories. It is proved that the estimated biases  converge  towards the true ones, if and only if, the scene admits no cylindrical symmetry.  The observer design  can be adapted to realistic  sensors where brightness and depth data are only available on a subset of $\SSS^2$. Preliminary simulations with synthetic brightness  and depth images  (corrupted by noise around 10\%) indicate that such Lyapunov-based observers  should be  robust and  convergent  for much weaker regularity assumptions.
\end{abstract}

\pagestyle{myheadings}
\thispagestyle{plain}
\markboth{N.Zarrouati-Vissi\`{e}re, K.Beauchard and P.Rouchon}{Bias estimation}

\section{Introduction}
\subsection{Context and challenges}
The  problem of estimating the position and the orientation of a moving object such as a ground, an aerial or an underwater vehicle has been extensively studied since World War II. In the 1950's, expensive inertial measurement units (IMUs) were developed, as missile guidance and control required extremely accurate navigation data \cite{missile}. Less expensive, tactical grade IMUs enable dead-reckoning techniques over short time periods, but require position fixes provided by GPS \cite{navigationGPS}, or combination through data fusion of other sensors outputs \cite{Car_navigation_survey,vissiere-et-al-cdc07}. As to recent low-cost IMUs using MEMS (Microelectromechanical systems) technologies, the cumulated error due to the bias of gyroscopes integrated over long time periods induces drift in orientation.  This drift can be managed; but from accelerometers, only high frequency output (dynamics) can be relied on. As odometers and velocimeters (e.g. Doppler radar \cite{navigationDoppler}, Pitot tube, electromagnetic (EM) log sensor), are commonly available technologies in vehicles, mass market applications can combine their linear velocity outputs with angular velocity from low-quality IMUs. Unfortunately, Pitot probes and EM log sensors are known to only provide airspeed and speed-through-the-water (STW) instead of speed-over-ground (SOG). We intend to study the situation, where linear and angular velocity are provided up to a slowly varying bias (the wind or the ocean current), which can be illustrated by \cite{bristeau2011ARdrone}.

Our approach leans on vision techniques: the field of dynamic vision mainly focuses on the estimation of motion of a camera and structure of a scene from a sequence of images \cite{pollefeys2004visual,Faugeras}. It usually tracks feature points between images and simultaneously recovers their three-dimensional (3D) position and the ego-motion of the camera through extended Kalman filtering (simultaneous localization and motion, SLAM) \cite{soatto1997} or non-linear observers \cite{Gupta06,AbdursulIG04,Astolfi10}.
Two difficulties systematically arise in those methods: as monocular systems can only estimate translation up to a scale factor, an additional output is required \cite{hamelhomography}; perspective systems induce nonlinearities to the system dynamics, which forces to study other geometrical formulations (e.g the essential space \cite{Motion_dynamic_vision}, the Pl\"{u}cker coordinates \cite{mahony2005image}).

The Kinect device has been a huge outbreak in the robotics and vision communities (\cite{kinectfusion}) as it provides depth measurements registered at each pixel of a RGB image, at a relatively low cost.  It enables the simultaneous exploitation of image and depth as dense data rather than sparse features, which to our knowledge have rarely been attempted. Such method, well-known as optical flow \cite{HS81} when it deals with image data, can by extension be considered as geometrical flow when it exploits depth images.

The contributions of this paper can be summarized as follows. We propose an original method to  estimate constant biases on angular and translational velocities and, at the same time, filter the brightness and depth  images. This method relies on
an $SO(3)$-invariant partial differential system~\cite{Bonnabel-Rouchon2009,ACC12_SO3} coupling, for a static scene,   the brightness field perceived by a spherical camera,  the depth field  and the angular/translational velocities.  Observability of this problem is studied. A precise geometric characterization of scenes that prevent observability of the biases is given in Theorem~\ref{thm:IS}: the biases are observable, if, and only if the scene does not admit a cylindrical symmetry.   The observer design is based on a Lyapunov functional and is also $SO(3)$-invariant. It yields an  integro/partial differential system for the estimated fields and biases (see~\eqref{eq:observer} for   $\SSS^2$   and see~\eqref{eq:obs_reduit} for  localization  on a spherical cap).   Asymptotic convergence is investigated under $C^1$ regularity assumptions on  motions and scene (Theorem~\ref{thm:CV}). The functional analysis relies on the adaptation of usual arguments  to our non-linear Partial Differential Equation  with non-local terms. The observer design and its  convergence  analysis  fully exploit   invariant differential calculus \cite{boothby-book} on the Riemannian sphere $\SSS^2$ where the independent variables labeling the image pixels live.
Simulations (figures~\ref{fig:biaistransbruit} and~\ref{fig:biaisrotbruit}) show that such observers could be used with noisy fields and  even without $C^0$ regularity of the entire scene.

Section~\ref{sec:SO3model} is devoted to  the $SO(3)$-invariant model,  the  partial differential system   coupling  brightness and  depth fields and the $C^1$ regularity and geometric assumptions used for the observer convergence analysis. In section~\ref{sec:observability} we prove that the biases are observable if, and only if the scene admits no cylindrical symmetry.  In section \ref{sec:asymp_obs},  the Lyapunov-based nonlinear observer is introduced and its convergence is investigated.  In section \ref{sec:implem}, we  explain how to adapt this observer to a realistic pinhole camera model with restricted fields of view and we present simulations  illustrating convergence and  robustness to noise of 10\%.

\subsection{Notations}
\begin{remunerate}
\item $C^k_b([0,+\infty),X)$ is the set of functions whose time derivatives are bounded on $X$ uniformly with respect to $t \in [0,+\infty)$ up to order $k$ .
\item The Euclidean scalar product of two vectors $\bold{a}$ and $\bold{b}$ in $\R^3$ is denoted by $\bold{a}\cdot \bold{b}$ and their wedge product by $\bold{a}\times \bold{b}$.

\item If $\SSS^2\ni \etab\mapsto f(\etab)\in\R$ is a scalar $C^1$ field on $\SSS^2$,
then $\nabla f$ denotes  the gradient of $f$ on the Riemannian unitary sphere $\SSS^2$ of $\mathbb{R}^3$:
for each vector $\etab \in \SSS^2$, $\nabla f(\etab)$ is a vector of  $\mathbb{R}^3$, orthogonal to $\etab$ and thus tangent to $\SSS^2$ at $\etab$.
\item If $\SSS^2\ni \etab\mapsto \bold{f}(\etab)\in\R^3$ is a $C^1$ vector field on $\SSS^2$ ($\forall \etab\in\SSS^2$, $\bold{f}(\etab)\cdot\etab \equiv 0$),
then $\nabla\cdot \bold{f}$ denotes  the divergence of $\bold{f}$ on the Riemannian sphere $\SSS^2$.
\end{remunerate}
\section{The $SO(3)$-invariant model }
\label{sec:SO3model}

\subsection{Modelling and regularity assumptions}
\label{subsec:modelling}

The model is based on geometric  assumptions introduced  in  \cite{Bonnabel-Rouchon2009,ACC12_SO3}. They are recalled in this sub-section. The camera is  spherical. Its motion is given through the linear and angular velocities $\vit(t)\in\R^3$ and ${\om} (t)\in\R^3$ expressed in a reference frame attached to the camera: the camera frame. More precisely, the position of the optical center  in the reference frame $\cal R$ is denoted by $C(t)\in\R^3$. Orientation versus $\cal R$  is given by the unitary quaternion $q(t)$: any vector ${\varsigma}$ in the camera frame corresponds to the vector $q{\varsigma} q^*$ in the reference frame $\cal R$ using the identification of  vectors of $\R^3$ as imaginary quaternions. $q^*$ denotes the conjugate of $q$ ($q^* q = qq^*=1$). We have thus
\begin{equation}\label{eq:chgtrepere}
\dot q = \frac{1}{2} q {\om} \text{ and } \dot C = q {\vit} q^*.
\end{equation}
A pixel is labeled by  a vector $\etab\in\SSS^2$  in the camera frame and  receives the brightness $y(t,\etab)$. Thus at each time $t$, the image produced by the camera is described by the  scalar field $\mathbb S^2\ni \etab \mapsto y(t,\etab)\in\R$.

The scene is modeled as a closed, $C^1$ and convex  surface $\Sigma$  of $\R^3$,  diffeomorphic to  $\mathbb  S^2$.   The camera is inside the domain $\Omega\subset \R^3$ delimited by $\Sigma=\partial\Omega$.
To a point $M\in\Sigma$ corresponds  one and only one camera pixel. At each time $t$, there is a  bijection between the  position of the pixel given by $\etab\in\SSS^2$ and the point $M\in\Sigma$. Since the point $M$ is labeled by $s\in\SSS^2$, this means that for each $t$, exist two mappings  $\SSS^2\ni s \mapsto \etab=\phi(t,s)\in\SSS^2$ and  $\SSS^2\ni \etab \mapsto s=\psi(t,\etab)\in\SSS^2$ with
$\phi(t,\psi(t,\etab))\equiv \etab$ and $\psi(t,\phi(t,s))\equiv s$, for all $\etab,s\in\SSS^2$. To summarize we have:
\begin{equation} \label{diffeo}
\etab= \phi(t,s)= q(t)^* \frac{\overrightarrow{C(t)M(s)}}{\| \overrightarrow{C(t)M(s)}\|} q(t)
\quad \text{and}\quad
s=\psi(t,\etab)
\end{equation}
where $\psi(t,.)$ and $\phi(t,.)$ are diffeomorphisms of $\SSS^2$ for every $t\geq 0$.

\begin{figure}[thpb]
      \centering
      \includegraphics[scale=0.4]{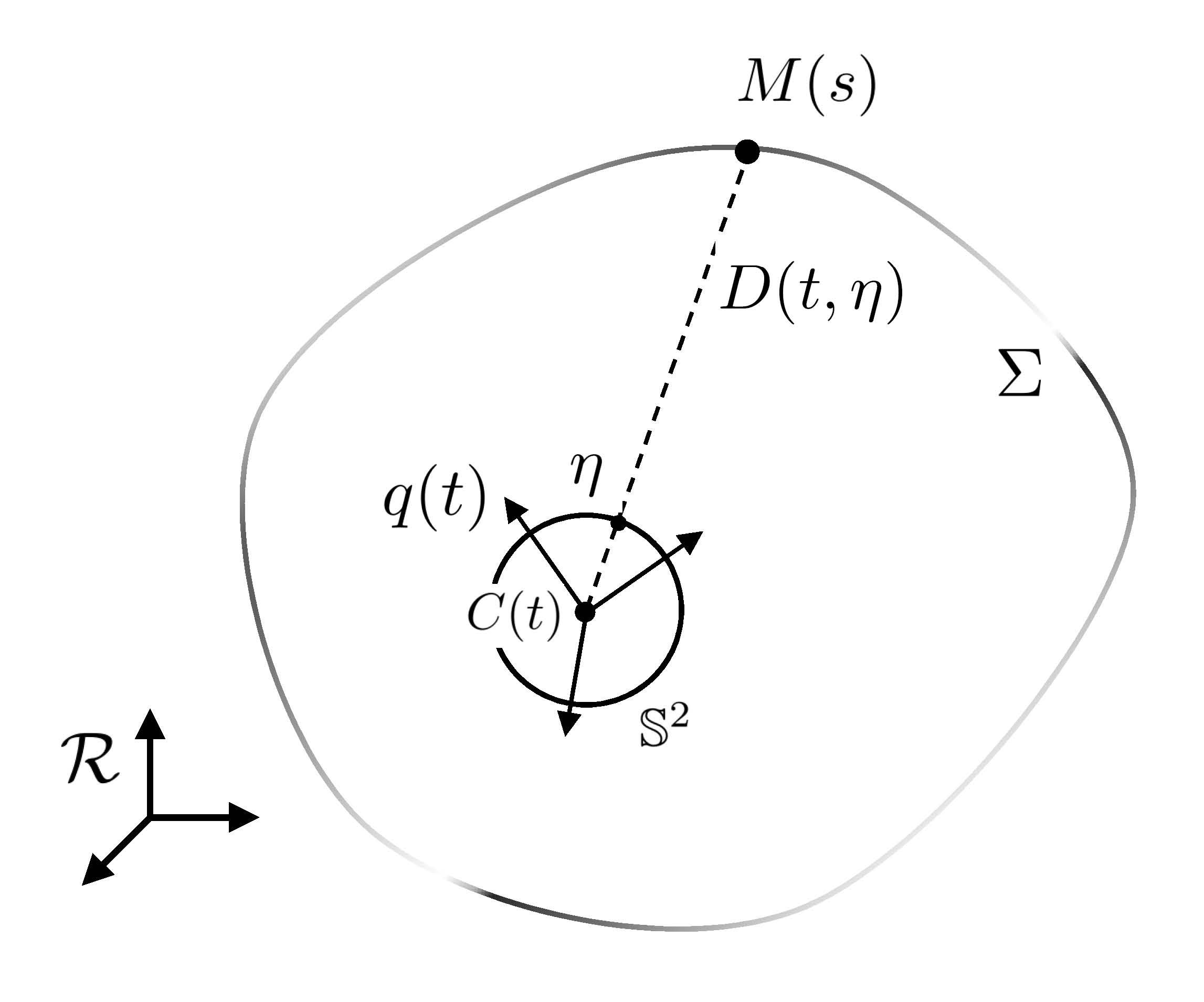}
      \caption{Model and notations of a spherical camera in a static scene described by a fixed, convex, closed and Lambertian  surface $\Sigma$.  \cite{Bonnabel-Rouchon2009,ACC12_SO3}.}
      \label{fig:notations}
\end{figure}

The intensity  of light emitted by a point  $M(s)\in\Sigma$ does not depend on the direction of emission ($\Sigma$ is a Lambertian surface) and is independent of $t$ (the scene is static). This means that $y(t,\etab)$ depends only on $s$:
there exists a function $y_\Sigma(s)$ such that
\begin{equation} \label{y_explicit}
y( t , \etab )= y_\Sigma ( \psi(t,\etab) ).
\end{equation}
We denote by $D_\Sigma(t,s)$ the distance between $C(t)$ and $M(s)$, $s\in\SSS^2$.
 Thus   the distance $D(t,\etab)$  between the optical center and the object seen in the direction $\etab\in\SSS^2$ is given by
\begin{equation} \label{D_explicit}
D(t,\etab)=D_\Sigma(t,\psi(t,\etab))=\| \overrightarrow{C(t)M(\psi(t,\etab))}\|.
\end{equation}
Fig.\ref{fig:notations} illustrates the model and the notations.

We assume
\begin{remunerate}
\item $\vit$ and $\om$ are in $C^1_b([0,+\infty),\R^3)$;  $C(t)$ stays in a fixed compact subset $K$ of $\Omega$ for all $t>0$,
\item $y_\Sigma$ is a $C^1$ function of $s \in \SSS^2$;  the surface $\Sigma$ is $C^1$.
\end{remunerate}
Thanks to the first  assumption and  \eqref{eq:chgtrepere},  $C(t)$ and $q(t)$ belong to $C^2_b([0,+\infty))$.
Thanks to the second assumption  and \eqref{diffeo}, $\psi(t,\etab)$ and $\phi(t,s)$ belong to $C^1_b([0,+\infty)\times\SSS^2,\SSS^2)$. We deduce also  from
\eqref{y_explicit}, \eqref{D_explicit} and the first assumption that
\begin{equation}\label{yDC1b}
\begin{array}{c}
y, D \in C^1_b([0,+\infty)\times\SSS^2,\R)\, \text{ and }
\\
D(t,\etab) \geqslant D_*>0\,, \forall (t,\etab) \in [0,+\infty)\times\SSS^2\,.
\end{array}
\end{equation}

\subsection{Statement of the bias estimation problem}
\label{subsec:kinem}
Under the above  assumptions, the functions $y$ and $D$ obey to (\cite{Bonnabel-Rouchon2009,ACC12_SO3})
\begin{align}
\partial_t y &=- \nabla y \cdot(\etab \times (\om +\frac{1}{D}\etab\times \vit)),
\label{eq:flot}
\\
\partial_t D &=- \nabla D \cdot(\etab \times (\om +\frac{1}{D}\etab\times \vit)) -\vit\cdot\eta,
\label{eq:prof}
\end{align}
where $\partial_t y$ and $\partial_t D$ stand for  partial derivatives of $y$ and $D$ with respect to $t$.
Equations \eqref{eq:flot} and \eqref{eq:prof} are $SO(3)$-invariant in the following sense: they remain unchanged by  any rotation described by the quaternion $\sigma$ and changing $(\etab,\om,\vit)$ to
$(\sigma \etab \sigma^*,\sigma \om \sigma^*,\sigma \vit \sigma^* )$.

The camera motion characterized by $(\vit,\om)$ is not known  precisely. We assume here  that sensor data provide $\vm(t)$ and $\wm(t)$ differing from true velocities by measurement biases, defined as constant errors:
\begin{equation}
\vm(t)=\vit(t)+\pv, \quad
\wm(t)=\om(t)+ \pw ,
\label{eq:vit_biais}
\end{equation}
where  $\pv, \pw \in\mathbb{R}^3$ are constant:
\begin{equation}
\partial_t \pv =0, \quad
\partial_t \pw =0.
\label{eq:biais_const}
\end{equation}

Thus, the state equations of the system are
\begin{equation}
 \left\{
    \begin{array}{l}
			\partial_t y=-\nabla y \cdot \etab \times(\wm-\pw+\frac{1}{D} \etab \times (\vm-\pv)),\\
      \partial_t D=-\nabla D \cdot \etab \times(\wm-\pw+\frac{1}{D} \etab \times (\vm-\pv))- \etab \cdot (\vm-\pv),\\
      \partial_t \pv=0,\\
      \partial_t \pw=0.
\end{array}
\right.
\label{eq:state_equations}
\end{equation}

 Here $(\vm,\wm)$ are considered as known inputs and  $(y,D)$ as  the measured outputs. The state  is $(y,D,\pv,\pw)^T$. The goal  is to estimate  in real-time the parameters $(\pv,\pw)$  from the known signals corresponding to sensors measuring  $y$, $D$, $\vm$ and $\wm$.

\section{Observability of  system~\eqref{eq:state_equations} with $(y,D)$ as measured output}
\label{sec:observability}

\subsection{Characterization of the observability}

To prove the observability, we consider two state-trajectories
$$
t\mapsto (y^1(t,\cdot),D^1(t,\cdot),\pv^1(t),\pw^1(t))^T \text{ and } t\mapsto (y^2(t,\cdot),D^2(t,\cdot),\pv^2(t),\pw^2(t))
$$
and assume that they satisfy the same state equations \eqref{eq:state_equations} with the same inputs $t\mapsto (\vm(t),\wm(t))$ and outputs $t\mapsto (y(t,\cdot),D(t,\cdot))$.
It immediately follows that $y^1=y^2=y$, $D^1=D^2=D$, $\pv^1=\pv^2+\pv$ and $\pw^1=\pw^2+\pw$ with $\partial_t \pv=0$ and $\partial_t \pw=0$.
Then $t\mapsto (y(t,\cdot),D(t,\cdot),\pv(t),\pw(t))^T$ satisfy:
\begin{equation}
 \left\{
    \begin{array}{l}
\nabla y \cdot \etab \times(\pw+\frac{1}{D} \etab \times \pv)=0,\\
\nabla D \cdot \etab \times(\pw+\frac{1}{D} \etab \times \pv)+\etab\cdot \pv=0,\\
\partial_t \pv=0,\\
\partial_t \pw=0.
\end{array}
\right.
\label{eq:invariant_space_extended}
\end{equation}
Observability means that \eqref{eq:invariant_space_extended}  implies that   $(\pw,\pv)=0$.

Consider a solution $(y,D,\pw,\pv)$ of~\eqref{eq:invariant_space_extended}. Take some fixed value $\bar t$ of $t$ and set $\bar y(\cdot)=y(\bar t,\cdot)$ with
$\bar D(\cdot)=D(\bar t,\cdot)$. Then $(\bar y,\bar D)$ is a stationary solution of~(\ref{eq:flot},\ref{eq:prof}) when  $(\om,\vit)$ is constant and equal to $(\pw,\pv)$.  Characterizing a scene, i.e., the surface $\Sigma$ and its painting  $y_{\Sigma}$,   that admits   a stationary solution of~\eqref{eq:flot} and~\eqref{eq:prof} with non zero  constant values of $(\om,\vit)$, i.e.  that is compatible with   an effective motion  producing stationary   images $y$ and $D$, is not obvious.

 These considerations motivate the following terminology: for  a given scene $(\Sigma,y_{\Sigma})$,  a constant  motion $(\om,\vit)=(\pw,\pv)$,  such that  a stationary solution $(y,D)$  of~(\ref{eq:flot},\ref{eq:prof}) exists, i.e.,
\begin{equation}
 \left\{
    \begin{array}{l}
\nabla D \cdot \etab \times(\pw+\frac{1}{D} \etab \times \pv)+\etab\cdot \pv=0,\\
\nabla y \cdot \etab \times(\pw+\frac{1}{D} \etab \times \pv)=0,\\
\end{array}
\right.
\label{eq:invariant_space}
\end{equation}
is called a stationary  motion inside   $(\Sigma,y_\Sigma)$.  For any scene  $(\Sigma,y_\Sigma)$, the  null motion $\pv=\pw=0$ is always stationary. Apart from this trivial stationary  motion, only very specific scenes (see Theorem~\ref{thm:IS})  admit non trivial stationary   motions $(\pv,\pw)$.
 Thus observability of~\eqref{eq:state_equations}  with output $(y,D)$   is equivalent to the following  statement: the stationary motions $(\pv,\pw)$ inside the  scene $(\Sigma,y_{\Sigma})$ are all  reduced to the trivial one $(0,0)$.

\subsection{Scenes   with non trivial stationary motion}

The following theorem characterizes the scenes admitting only trivial stationary motion, i.e., such that a stationary solution of~(\ref{eq:flot},\ref{eq:prof}) exists for some  $(\om,\vit)\neq 0$ constant.
\begin{theorem} \label{thm:IS}
Consider the assumptions given in subsection~\ref{subsec:modelling}. The scene  $(\Sigma,y_{\Sigma})$ admits a non trivial stationary motion, if and only if, it admits  a rotation axis.
\end{theorem}
\\
{\em Proof of Theorem~\ref{thm:IS}:}
Let us consider a non-zero motion $(\pv,\pw)^T$. Under the assumption that the equations \eqref{eq:invariant_space} have a solution $(y,D) \in C^1(\SSS^2,\R^2)$ we prove that $(\Sigma,y_{\Sigma})$ admits a rotation axis.
Note that
\begin{align*}
&\Sigma=\left\{D(\etab)\etab;\etab\in\SSS^2\right\},\\
&y_{\Sigma}(s)=y(t,\etab) \text{ with } \etab=\phi(t,s) \text{(see Section \ref{subsec:modelling})}.
\end{align*}
Let us eliminate two trivial cases:
\begin{romannum}
\item $\pv=0, \pw\neq0$: the two first equations of \eqref{eq:invariant_space} reduce to
\begin{equation}
 \left\{
    \begin{array}{l}
\nabla D \cdot (\etab \times\pw)=0,\\
\nabla y \cdot (\etab \times\pw)=0,\\
\end{array}
\right.
\end{equation}
which characterizes a cylindrical symmetry for $y$ and $D$, with same rotation axis defined by the direction of $\pw$;
\item $\pv\neq0, \pw=0$: let us replace $\etab$ by $\pv/\|\pv\|$ in the first equation of \eqref{eq:invariant_space}: this yields $\|\pv\|=0$, which is absurd.
\end{romannum}
From now on, we consider only the case where $\pv$ and $\pw$ are both non-zero.
Our strategy consists in finding a fixed vector $\bu$ such that a translation of the origin of the camera frame in the direction $\bu$ transforms the functions $D(\etab)$ and $y(\etab)$ in new functions $\tD$ and $\ty$ on the sphere $\SSS^2$ and those functions are invariant by rotation about the axis $\R\pw$. This requires several technical steps.

\textbf{Step 1: Let us prove that $\pv$ and $\pw$ are orthogonal.} Starting from the first equation of~\eqref{eq:invariant_space}, let us multiply it by $D^n(\etab\cdot\pw)$, where the choice of $n \in \NN$ will follow, and integrate the result on $\SSS^2$:
\begin{equation}
I=\int_{\SSS^2}\Big( D^n\nabla D \cdot \etab \times(\pw+ \frac{1}{D} \etab \times \pv)+D^n\eta\cdot \pv \Big) (\etab\cdot\pw)d\sigma_\eta=0,
\end{equation}
We have
\begin{multline*}
I=\int_{\SSS^2} \Big[\frac{1}{n+1}\nabla D^{n+1} \cdot (\etab \times\pw)(\etab\cdot\pw)
+\frac{1}{n} \nabla D^n \cdot (\etab \times(\etab\times\pv))(\etab\cdot\pw)\\
+D^n(\etab\cdot \pv)(\etab\cdot\pw)\Big]d\sigma_\eta.
\end{multline*}
The two first terms are integrated by parts:
\begin{multline*}
I=\int_{\SSS^2}\Big[-\frac{1}{n+1} D^{n+1} \nabla\cdot ((\etab \times\pw)(\etab\cdot\pw))
-\frac{1}{n}  D^n \nabla\cdot ((\etab \times(\etab\times\pv))(\etab\cdot\pw))\\
+D^n(\etab\cdot \pv)(\etab\cdot\pw)\Big]d\sigma_\eta.
\end{multline*}
The derivatives are developed:
\begin{multline*}
I=\int_{\SSS^2}\Big[ -\frac{1}{n+1} D^{n+1}  (\nabla\cdot(\etab \times\pw)(\etab\cdot\pw)+(\etab \times\pw)\cdot\nabla(\etab\cdot\pw))\\
 -\frac{1}{n}  D^n \Big( \nabla\cdot (\etab \times(\etab\times\pv))(\etab\cdot\pw)+(\etab \times(\etab\times\pv))\cdot\nabla(\etab\cdot\pw) \Big)\\
 +D^n(\etab\cdot \pv)(\etab\cdot\pw)\Big]d\sigma_\eta.
\end{multline*}
Now, let us recall the following basic formulae of differential geometry on $\SSS^2$ (detailed and proven in appendix), where $\bold{P}$ is a constant vector:
\begin{align}
\nabla(\etab\cdot \bold{P})&=-\etab\times(\etab\times \bold{P}),
\label{formule:gradient_scalaire}
\\
\Delta(\etab\cdot \bold{P})&=-2\etab\cdot \bold{P},
\label{formule:laplacien}
\\
\nabla\cdot(\etab\times \bold{P})&=0.
\label{formule:divergence}
\end{align}
Using the formulae \eqref{formule:divergence} and \eqref{formule:gradient_scalaire} in the first term, and \eqref{formule:laplacien} and \eqref{formule:gradient_scalaire} in the second term yields:
\begin{multline*}
I=\int_{\SSS^2}\Big[\frac{1}{n+1} D^{n+1} (\etab \times\pw)\cdot(\etab \times(\etab\times\pw))\\
-\frac{1}{n}  D^n ((2\etab \cdot\pv)(\etab\cdot\pw)-(\etab \times(\etab\times\pv))\cdot(\etab \times(\etab\times\pw)))\\
+D^n(\etab\cdot \pv)(\etab\cdot\pw)\Big]d\sigma_\eta.
\end{multline*}
The first term is obviously zero. Then
$$
(\etab \times(\etab\times\pv))\cdot(\etab \times(\etab\times\pw))=-(\etab\cdot\pv)(\etab\cdot\pw)+\pv\cdot\pw.
$$
Thus
$$
I=\int_{\SSS^2}\Big[-\frac{1}{n}  D^n ((3\etab \cdot\pv)(\etab\cdot\pw)-\pv\cdot\pw)+D^n(\etab\cdot \pv)(\etab\cdot\pw)\Big]d\sigma_\eta.
$$
Choosing $n=3$ yields
\begin{equation}
(\pv\cdot\pw)\int_{\SSS^2} D^3 d\sigma_\eta=0.
\end{equation}
Since $D>0$, we deduce that $\pv\cdot\pw=0$. This ends the proof of Step 1.

As $\pv$ and $\pw$ are orthogonal and both non-zeros, one can introduce a new vector, $l\bu$, where $l>0$ is a length and $\bu$ is a unit vector, by: $\pv=\pw\times l\bu$. Let us also define $\bv$ and $\bw$ the unit vectors such that $\pv=\|\pv\| \bv$ and $\pw=\|\pw\| \bw$. Then, $(\bu,\bv,\bw)$ is an orthonormal and direct frame of $\R^3$.
From now on, \eqref{eq:invariant_space}, multiplied by $lD/\|\pv\|$ writes:
\begin{equation}\label{eq:ISD_orthogonal}
\left\{
    \begin{array}{l}
\nabla D \cdot ((D\etab+l\bu)\times \bw)+lD\eta \cdot \bv=0,\\
\nabla y \cdot ((D\etab+l\bu)\times \bw)=0.
\end{array}
\right.
\end{equation}
\\

\textbf{Step 2: Let us prove that $D(\etab)\etab+l\bu\neq0, \forall\etab\in\SSS^2$.}
It is sufficient to prove that $D(-\bu)\neq l$.
Let $\etab_\epsilon:=\frac{-\bu+\epsilon\bv}{\sqrt{1+\epsilon^2}}=-\bu+\epsilon \bv + o(\epsilon)$.
We have
$$\begin{array}{ll}
D(\etab_\epsilon)\etab_\epsilon
& = \left(D(-\bu)+\epsilon\left.\nabla D\right|_{-\bu}\cdot\bv+o(\epsilon)\right)\left( -\bu+\epsilon \bv + o(\epsilon)  \right)\\
& = -D(-\bu)\bu+\epsilon\left(D(-\bu)\bv-\left(\left.\nabla D\right|_{-\bu}\cdot\bv\right)\bu\right)+o(\epsilon)
\end{array}$$
Thus, \eqref{eq:ISD_orthogonal} evaluated at $\etab_\epsilon$ yields
\begin{align*}
\left.\nabla D\right|_{\etab_\epsilon} \cdot \left[\left[(l-D(-\bu))\bu +\epsilon\left(D(-\bu)\bv-\left(\left.\nabla D\right|_{-\bu}\cdot\bv\right)\bu\right)  \right]\times \bw\right]+\epsilon lD(-\bu) = o(\epsilon),
\end{align*}
or
\begin{align}\label{eq:gradient_epsilon_2}
\left.\nabla D\right|_{\etab_\epsilon} \cdot \left[(D(-\bu)-l)\bv +\epsilon\left(D(-\bu)\bu+\left(\left.\nabla D\right|_{-\bu}\cdot\bv\right)\bv\right)\right]+\epsilon lD(-\bu)
= o(\epsilon).
\end{align}
Equation \eqref{eq:ISD_orthogonal} evaluated at $\etab=-\bu$ and substracted to \eqref{eq:gradient_epsilon_2} yields:
\begin{multline}\label{eq:gradient_epsilon_3}
(D(-\bu)-l)\left(\left.\nabla D\right|_{\etab_\epsilon}-\left.\nabla D\right|_{-\bu} \right) \cdot \bv \\+\epsilon\left[D(-\bu)\left.\nabla D\right|_{\etab_\epsilon}\cdot\bu+\left(\left.\nabla D\right|_{-\bu}\cdot\bv\right)\left(\left.\nabla D\right|_{\etab_\epsilon}\cdot\bv\right)+ lD(-\bu)\right] = o(\epsilon).
\end{multline}
Let us define
$$\begin{array}{cccc}
F: & \R        & \rightarrow & \R \\
   & \epsilon  & \mapsto     &
D(-\bu)\left.\nabla D\right|_{\etab_\epsilon}\cdot\bu+\left(\left.\nabla D\right|_{-\bu}\cdot\bv\right)\left(\left.\nabla D\right|_{\etab_\epsilon}\cdot\bv\right)+ lD(-\bu).
\end{array}$$

Assuming that $(D(-\bu)-l)=0$ implies that $F(\epsilon)=o(1) $ as  $\epsilon\rightarrow 0$. Since $F$ is a $C^0$ function of $\epsilon$, it implies $F(0)=0$, i.e.
$$\left(\left.\nabla D\right|_{-\bu}\cdot\bv\right)^2+lD(-\bu)=0$$ which is absurd since $l>0$. This concludes the Step 2.
\\

\textbf{Step 3: Let us prove that $f(\etab):=\|D(\etab)\etab+l\bu\|$ (respectively $y(\etab)$) satisfies $df(\etab)\cdot\delta\etab=0$ (respectively $dy(\etab)\cdot\delta\etab=0$) where $\delta\etab:=(D(\etab)\etab+l\bu)\times\bw+l(\etab\cdot\bv)\etab$,  $\forall \etab\in \SSS^2$.}
We have
$$f(\etab)=g \circ k(\etab) \quad \text{ where } \quad g(\boldsymbol{x}):=\|\boldsymbol{x}\|  \quad \text{ and } \quad k(\etab):=D(\etab)\etab+l\bu.$$
Thanks to the Step 2, $f \in C^1(\SSS^2)$ and
$$df(\etab)\cdot\delta\etab = dg[k(\etab)]\cdot\Big( dk(\etab)\cdot\delta\etab \Big).$$
We have
$$\begin{array}{ll}
dk(\etab)\cdot\delta \etab
& = \left(\nabla D(\etab)\cdot\delta\etab\right)\etab + D(\etab) \delta\etab \\
& = \left( \nabla D(\etab)\cdot \left( [D(\etab)\etab+l\bu] \times \bw+ l(\etab\cdot \bv)\etab \right) \right) \etab \\
&\quad + D(\etab) \left( [D(\etab)\etab+l\bu] \times \bw + l(\etab\cdot \bv)\etab \right) \\
& = \left( \nabla D(\etab)\cdot \left( [D(\etab)\etab+l\bu] \times \bw \right) \right)\etab + D(\etab) \left( [D(\etab)\etab+l\bu] \times \bw + l(\etab\cdot \bv)\etab \right)
\end{array}$$
because $\nabla D(\etab) \perp \etab$. Thanks to~\eqref{eq:ISD_orthogonal}, we obtain
$$dk(\etab)\cdot\delta \etab  = D(\etab)  [D(\etab)\etab+l\bu] \times \bw.$$
Thus
$$df(\etab)\cdot\delta\etab  = \left\langle \frac{D(\etab)\etab+l\bu}{\|D(\etab)\etab+l\bu\|} , D(\etab)  [D(\etab)\etab+l\bu] \times \bw \right\rangle =0.$$
Similarly,
\begin{align*}
dy(\etab)\cdot\delta\etab &= \nabla y\cdot\delta\etab\\
&= \nabla y\cdot((D(\etab)\etab+l\bu)\times\bw+l(\etab\cdot\bv)\etab)\\
&= \nabla y\cdot((D(\etab)\etab+l\bu)\times\bw)\text{ because } \nabla y(\etab) \perp \etab\\
&=0 \text{ according to \eqref{eq:ISD_orthogonal}.}
\end{align*}
This ends the proof of Step 3.
\\

\textbf{Step 4: Let us prove that, when $\etab=\etab(\tau)$ satisfies $\partial_\tau \etab=[D(\etab)\etab+l\bu] \times \bw + l(\etab\cdot \bv)\etab$,
then $\widetilde{\etab}:=[D(\etab)\etab+l\bu]/\|[D(\etab)\etab+l\bu]\|$
is a periodic trajectory on $\SSS^2$,
describing a circle perpendicular to the axis $\mathbb{R}\bw$.}
Let $\etab_0 \in \SSS^2$. We consider the characteristic (see \cite{leveque1992numerical,serre1999systems} for a description of the method of characteristics)
$$\left\lbrace \begin{array}{l}
\frac{d\etab}{d\tau}(\tau) =  \Big( D[\etab(\tau)] \etab(\tau)+l\bu \Big) \times \bw + l \Big(\etab(\tau)\cdot \bv\Big) \etab(\tau), \\
\etab(0)=\etab_0.
\end{array}\right.$$
The right hand side is defined by a $C^1$ function of $\etab$ thus the Cauchy-Lipschitz theorem guarantees the existence
and uniqueness of a maximal solution. Moreover, no explosion is possible in finite time (the solution $\etab(\tau)$ lives on $\SSS^2$) thus
the maximal solution is defined for every $\tau \in \mathbb{R}$.
From Step 3, we deduce that the quantity
$\tau \mapsto \|D[\etab(\tau)]\etab(\tau)+l\bu\|$ is constant along the characteristics:
$$\|D[\etab(\tau)]\etab(\tau)+l\bu\| \equiv \widetilde{D} := \|D[\etab_0]\etab_0+l\bu\|, \quad \forall \tau \in \R.$$
Now let us consider
$$\widetilde{\etab}(\tau):=\frac{D[\etab(\tau)]\etab(\tau)+l\bu}{\|D[\etab(\tau)]\etab(\tau)+l\bu\|} = \frac{D[\etab(\tau)]\etab(\tau)+l\bu}{\widetilde{D}}.$$
Then, $\widetilde{\etab}(\tau)$ is a $C^1$-function and
$$\begin{array}{ll}
\widetilde{D} \frac{d\widetilde{\etab}}{d\tau} &= \left( \nabla D[\etab(\tau)]\cdot\frac{d\etab}{d\tau} \right) \etab(\tau) + D[\etab(\tau)] \frac{d\etab}{d\tau} \\
&= \left( \nabla D[\etab(\tau)]\cdot \Big( D[\etab(\tau)] \etab(\tau)+l\bu \Big) \times \bw \right) \etab(\tau) \\
&\quad + D[\etab(\tau)] \Big( \Big( D[\etab(\tau)] \etab(\tau)+l\bu \Big) \times \bw + l \Big(\etab(\tau)\cdot \bv \Big) \etab(\tau) \Big)\,.
\end{array}$$
Here, we have used the relation $\nabla D[\etab(\tau)]\cdot\etab(\tau)=0$ to get the last equality.
Thanks to~\eqref{eq:ISD_orthogonal}, we obtain
$$
\widetilde{D} \frac{d\widetilde{\etab}}{d\tau}
 = D[\etab(\tau)] \Big( D[\etab(\tau)] \etab(\tau)+l\bu \Big) \times \bw
= D[\etab(\tau)] \widetilde{D} \widetilde{\etab}(\tau) \times \bw.
$$
Thus, $\widetilde{\etab}$ turns about the axis $\R\bw$ as $\tau$ varies.
Note that the speed $D[\etab(\tau)]$ is bounded from below by a positive constant,
thus the trajectory of $\widetilde{\etab}$ is the (complete) circle on $\SSS^2$,
perpendicular to the axis $\R\bw$ and containing $\widetilde{\etab}_0:=[D(\etab_0)\etab_0+l\bu]/\widetilde{D}$.
\\
\\

\begin{figure} \label{fig:translationdomain}
      \centering
     \includegraphics[scale=0.5]{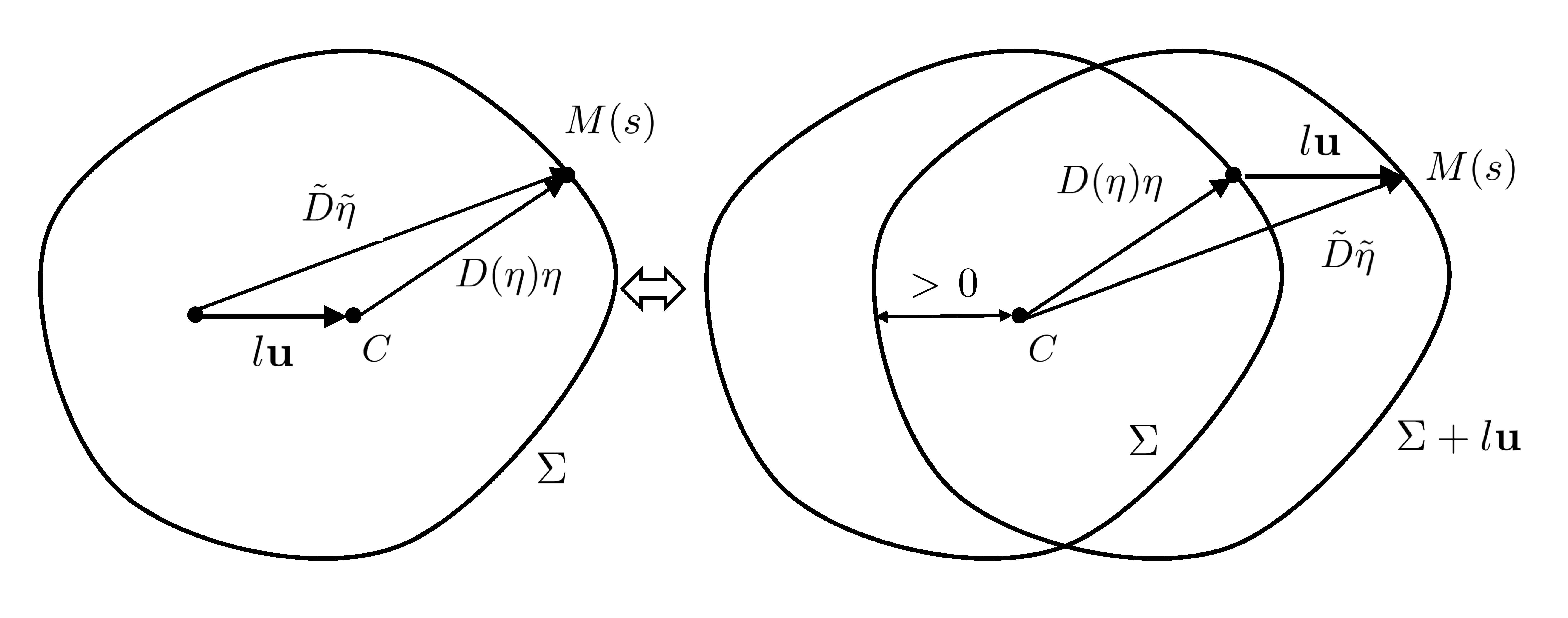}
      \caption{The vectors $D(\etab)\etab$ taken from the camera center $C$ and $\tD\widetilde{\etab}$ taken from the camera center translated by $-l\bu$ arrive at the same point $M(s)$. It is equivalent to see this as a translation of the domain $\Sigma$ by $l\bu$.}
\end{figure}

\textbf{Step 5: Let us prove that the map}
$$\begin{array}{ccc}
\SSS^2 & \rightarrow & \SSS^2 \\
\etab_0       & \mapsto     & \widetilde{\etab}_0:=\frac{D(\etab_0)\etab_0+l\bu}{\|[D(\etab_0)\etab_0+l\bu]\|}
\end{array}$$
\textbf{is surjective.}
This step leans on the statement that when the convex set defined by the surface $\Sigma$ is translated by $l\bu$, the center of the camera frame stays inside the convex set ($\Sigma+l\bu$). This statement is ensured by the fact that $D(-\bu) > l$ (or equivalently, that $\tD(-\bu) > 0$ , see Figure \ref{fig:translationdomain}).

Thus, for any $\widetilde{\etab}_0 \in \SSS^2$, $\exists M$ a point of the surface $\Sigma+l\bu$ and $\exists \lambda >0$ such that $\lambda \widetilde{\etab}_0=\overrightarrow{CM}=D(\etab_0)\etab_0+l\bu$ where $\etab_0 \in \SSS^2$ (by definition of $\Sigma+l\bu$). Then $\lambda=\|D(\etab_0)\etab_0+l\bu\|$ since $\widetilde{\etab}_0$ is a unit vector, and we have shown that $\exists \etab_0 \in \SSS^2$ such that $\widetilde{\etab}_0=\frac{D(\etab_0)\etab_0+l\bu}{\|D(\etab_0)\etab_0+l\bu\|}$. Thus the map $\etab_0\mapsto\widetilde{\etab}_0$ is surjective.

Now, let us prove that $D(-\bu) > l$. Let us assume that $D(-\bu) \leqslant l$. Then, the function
$\etab \mapsto D(\etab)+l\bu\cdot\etab$ is positive at $\etab=\bu$ and non positive at $\etab=-\bu$,
thus (intermediate values theorem) there exists $\etab^*$ such that $\etab^*\cdot \bw=0$ and
$D(\etab^*)+l\bu\cdot\etab^*=0$. Then, the vector $D(\etab^*) \etab^* + l\bu$ is orthogonal to
$\bw$ (because $\bu$ and $\etab^*$ are) and orthogonal to $\etab^*$, thus
$$D(\etab^*) \etab^* + l\bu = \left( \Big( D(\etab^*) \etab^* + l\bu \Big)\cdot (\bw \times \etab^*) \right) \bw \times \etab^*.$$
Note that $\bw \times \etab^*$ is unitary.
Moreover,
$$\begin{array}{ll}
\Big( D(\etab^*) \etab^* + l\bu \Big)\cdot (\bw\times \etab^*)
   & =l \bu.( \bw \times \etab^*) \text{ because } \etab^* \perp \bw \times \etab^*
\\ & =l \etab^*\cdot(\bu \times \bw)
\\ & =-l \etab^*\cdot \bv
\end{array}$$
thus
$$D(\etab^*) \etab^* + l\bu = l (\etab^*\cdot\bv) \etab^* \times  \bw .$$
Incorporating this relation in (\ref{eq:ISD_orthogonal}), we get
$$l (\etab^*\cdot\bv) \nabla D(\etab^*)\cdot \Big(  (\etab^* \times  \bw) \times \bw \Big) + D(\etab^*)l \etab^*\cdot\bv =0.$$
Moreover,
$$ \Big(  \etab^* \times  \bw \Big) \times \bw
=\etab^* - (\bw\cdot\etab^*)\bw = \etab^*$$
is orthogonal to $\nabla D(\etab^*)$ thus
$$l D(\etab^*) \etab^* \cdot\bv =0.$$
Therefore, $\etab^*\cdot\bw=\etab^*\cdot\bv=0$, i.e. $\etab^*=\pm \bu$, which is impossible because $D(\bu)>0$ and $D(-\bu)\neq l$ (see Step 2). This proves that $D(-\bu)>l$ and ends the Step 5.
\\

\textbf{Step 6: Conclusion.}
The origin of the camera frame is translated by the vector $-l\bu$. We denote by $\widetilde{\etab} \in \SSS^2$ the space variable in this new system of coordinates. Then, to any pair $(\etab,D(\etab))$ in the first system of coordinates corresponds a unique pair $(\widetilde{\etab},\tD)$  in the new system of coordinates:
$$\widetilde{\etab}:=\frac{D(\etab)\etab+l\bu}{\| D(\etab)\etab+l\bu \|}, \qquad
\widetilde{D} := \| D(\etab)\etab+l\bu \|,$$
and they describe the same point $M(s)$ on $\Sigma$:
$$M(s)=D(\etab)\etab=-l\bu + \widetilde{D} \widetilde{\etab}.$$
In addition, $\ty$ is defined so that for a given point $M(s)$ on $\Sigma$, the light intensity is the same in both systems of coordinates:
\begin{align*}
y_\Sigma(s)&=y(t,\etab) \text{ where }\etab=\phi(t,s)\text{ in the first system of coordinates}\\
&=\ty(t,\widetilde{\etab}) \text{ where } \widetilde{\etab}=\frac{D(\etab)\etab+l\bu}{\|D(\etab)\etab+l\bu\|}\text{ in the new system of coordinates.}
\end{align*}
Steps 3 and 4 of this proof show that the quantities $\tD$ and $\ty$ are constant when $\widetilde{\etab}$ covers circles on $\SSS^2$ orthogonal to the axis $\R\bw$. Moreover, Step 5 ensures that $\widetilde{\etab}$ covers any such circle.
Therefore, $\Sigma=\{\tD(\widetilde{\etab})\widetilde{\etab};\widetilde{\etab}\in\SSS^2\}$ (in the new system of coordinates) and $y_\Sigma(s)$
have rotation axis, with direction $\bw$. \hfill $\Box$

\section{The asymptotic observer}
\label{sec:asymp_obs}

\subsection{A Lyapunov based  observer}
\label{ssec:obs1}
We propose the following observer for $\pw$ and $\pv$, inspired by the one proposed in \cite{Bonnabel-Rouchon2009}:
\begin{equation}
 \left\{
    \begin{array}{ll}
        \partial_t\yh=&- \nabla \yh \cdot(\etab \times (\wm-\pwh +\frac{1}{D}\etab\times (\vm-\pvh))) +k_y(y-\yh),\\
        \partial_t\Dh=& - \nabla \Dh \cdot(\etab \times (\wm-\pwh +\frac{1}{D}\etab\times (\vm-\pvh)))-(\vm-\pvh)\cdot\etab +k_D(D-\Dh),\\
        \frac{d \pwh}{dt}=&-k_{\omega}\text{\Large $\int$}_{\SSS^2} \Big( \lambda_y(\yh-y)(\nabla \yh\times\etab)+\lambda_D(\Dh-D)(\nabla \Dh\times\etab) \Big)d\sigma_\eta,\\
        \frac{d \pvh}{dt}=&-k_{v}\text{\Large $\int$}_{\SSS^2}\left( \lambda_D(\Dh-D)\etab +\lambda_y(\yh-y)(\frac{1}{D}\etab\times(\etab\times\nabla \yh))\right.\\
        &\qquad\quad\left.+\lambda_D(\Dh-D)(\frac{1}{D}\etab\times(\etab\times\nabla \Dh))  \right)d\sigma_\eta,
    \end{array}
\right.
\label{eq:observer}
\end{equation}
where $k_y$, $k_D$, $k_{\omega}$, $k_v$, $\lambda_y$, $\lambda_D >0$
are constant gains.
The choice of this observer is motivated by the following Lyapunov function candidate
\begin{equation} \label{eq:Lyapunov}
V:= \frac{1}{2} \int_{\SSS^2} \left( \lambda_y \ty^2 + \lambda_D \tD^2 \right) d\sigma_\eta +
\frac{\|\tpw\|^2}{2 k_{\omega}} + \frac{\|\tpv\|^2}{2 k_v}
\end{equation}
with
\begin{equation}\label{eq:def_errors}
\ty:=\yh-y, \quad\tD:=\Dh-D,\quad \tpw:=\pwh-\pw\, \text{ and }\, \tpv:=\pvh-\pv.
\end{equation}
Let us prove that $V$ decreases along the trajectories of (\ref{eq:observer}) for appropriate choices of $k_y$ and $k_D$..
Let us consider any scalar field $h(t,\etab)$ defined on $\SSS^2$ and its integral $\mathcal {H}$ on the unit sphere:
\begin{equation*}
\mathcal {H}=\int_{\SSS^2}h(t,\etab)^2d\sigma_\eta=\int_{\SSS^2}h(t,q^{*}\etab q)^2d\sigma_\eta
\end{equation*}
since $\etab \mapsto q^{*}\etab q$ is an isometry on $\SSS^2$. For $\etab$ a constant vector of the earth-fixed frame, $q^{*}\etab q$ is the same vector, expressed in the camera frame. Thus, $h(t,q^{*}\etab q)$ is the value of $h$ corresponding to a specific object $M(s)$ of the scene. This yields
\begin{align*}
&\dotex \mathcal{H}=\int_{\SSS^2}\dotex(h(t,q^{*}\etab q)^2)d\sigma_\eta\\
&=2 \int_{\SSS^2}h(t,q^{*}\etab q)\left.\dotex h(t,q^{*}\etab q)\right|_s d\sigma_\eta\\
&=2 \int_{\SSS^2}h(t,q^{*}\etab q)\left(\partial_t h+\nabla h \cdot(\etab \times \om)\right) d\sigma_\eta.
\end{align*}
One can apply this calculation rule to the scalar fields $\ty$ and $\tD$. Then, equations \eqref{eq:flot} and \eqref{eq:observer} yield
\begin{equation}
\partial_t\ty+\nabla \ty \cdot(\etab \times \om) = \nabla \yh \cdot(\etab \times (\tpw +\frac{1}{D}\etab\times \tpv))-\nabla\ty\cdot(\etab\times(\frac{1}{D}\etab\times \vit))-k_y\ty\,.
\label{eq:ytilde}
\end{equation}
Equations \eqref{eq:prof} and \eqref{eq:observer} yield
\begin{equation}
\partial_t \tD +\nabla \tD\cdot(\etab \times \om )=\nabla \Dh \cdot(\etab \times (\tpw +\frac{1}{D}\etab\times \tpv))+\etab\cdot\tpv
-\nabla\tD\cdot(\etab\times(\frac{1}{D}\etab\times \vit))-k_D\tD\,.
\label{eq:gamtilde}
\end{equation}
Now, using $\bold{a}\cdot(\bold{b}\times \bold{c})=(\bold{a}\times \bold{b})\cdot \bold{c}$ and the expressions of  $\frac{d \pwh}{dt}$ and $\frac{d \pvh}{dt}$, one gets
\begin{equation}
\dot{V}=-\int_{\SSS^2} \left(k_y\lambda_y\ty^2+k_D \lambda_D\tD^2\right)d\sigma_\eta
-\int_{\SSS^2}\left( \lambda_y\frac{\nabla\ty^2}{2}+\lambda_D\frac{\nabla\tD^2}{2} \right) \cdot W(t,\etab) d\sigma_\eta.
\label{eq:lyapunov_der_1}
\end{equation}
where $(\R,\SSS^2)\ni(t,\etab)\mapsto W(t,\etab)=\etab\times(\frac{1}{D}\etab\times v)\in\R^3$.
Integration by parts of the second term of this expression yields
\begin{equation}
\dot{V}=-\int_{\SSS^2} \left(\lambda_y\left(k_y-\frac{\nabla\cdot W(t,\etab)}{2}\right)\ty^2+ \lambda_D\left(k_D-\frac{\nabla\cdot W(t,\etab)}{2}\right)\tD^2\right)d\sigma_\eta\,.
\label{eq:lyapunov_der_2}
\end{equation}
Thus $V$ is a Lyapunov function for trajectories of (\ref{eq:observer}) when $k_y, k_D$ are large enough (see Proposition \ref{Prop:Lyapunov} below).
\\

The convergence analysis relies on an adaptation to this infinite dimensional setting of the Lasalle invariance principle.
Precisely, we prove that any adherence value of $(\ty,\tD,\tpw,\tpv)$ as $[t \rightarrow + \infty]$
satisfies geometric equations equivalent to \eqref{eq:invariant_space} (see Theorem \ref{thm:CV}).
This proof requires the following successive steps.
\begin{remunerate}
  \item The Lyapunov function $V$ decreases along the trajectories of~\eqref{eq:observer}, for appropriate choices of $k_y$ and $k_D$ (see Proposition \ref{Prop:Lyapunov}).
  \item The non-linear partial differential equations~\eqref{eq:observer} define a well posed Cauchy problem
for every $t\geq 0$ (see Proposition \ref{Prop:WP1}).
  \item The functions $\nabla\yh$ and $\nabla\Dh$ are uniformly bounded with respect to $(t,\etab) \in [0,+\infty) \times \SSS^2$ (see Proposition \ref{Prop:WP2}).
  \item The estimated state $(\yh,\Dh,\pwh,\pvh)$ depends continuously on $(y,D,\om,\vit)$ (see Proposition \ref{Prop:WP3}).
\end{remunerate}

\subsection{Decrease of the Lyapunov function}
In the next statement are given sufficient conditions for the choice of correction gains $k_y$ and $k_D$ to ensure the decrease of the Lyapunov function $V$.

\begin{proposition} \label{Prop:Lyapunov}
The Lyapunov function $V$ decreases along the trajectories of (\ref{eq:observer}) under the assumption that $k_y,k_D > \frac{L}{2}$
where $L$ is defined as
$$L:= \left\|
\frac{1}{D(t,\etab)}\left(2\etab\cdot \vit(t)-\frac{\nabla D(t,\etab)}{D(t,\etab)}\cdot (\etab\times(\etab\times \vit(t)))\right)
\right\|_{L^\infty((0,+\infty)\times\SSS^2)}.$$
\end{proposition}
{\em Proof of Proposition \ref{Prop:Lyapunov}:}
Note that $L$ is well defined thanks to (\ref{yDC1b}) and assumption $\vit \in C^1_b((0,+\infty),\mathbb{R}^3)$ (see subsection \ref{subsec:modelling}).
Moreover, $L=\|\nabla\cdot W(t,\etab)\|_{L^\infty((0,+\infty)\times\SSS^2)}$ (see formulae \eqref{formule:gradient_scalaire} and \eqref{formule:laplacien}).
From \eqref{eq:lyapunov_der_2}, we deduce that
\begin{equation}\label{eq:Lyapdecrease}
\dot{V}\leq-\int_{\SSS^2} \left(\lambda_y\left(k_y-\frac{L}{2}\right)\ty^2+ \lambda_D\left(k_D-\frac{L}{2}\right)\tD^2\right)d\sigma_\eta:=-f(t)
\end{equation}
which ends the proof of Proposition \ref{Prop:Lyapunov}. \hfill $\Box$
\\
From now on, we assume that $k_y$ and $k_D$ satisfy the assumption of Proposition \ref{Prop:Lyapunov}.

\subsection{Existence and uniqueness}
Existence and uniqueness of solutions of (\ref{eq:observer}) is given by the next statement.
For the simplicity of notations, from now on, we assume that $k_v=k_{\omega}=1$.

\begin{proposition} \label{Prop:WP1}
Let $R>0$. There exists a constant
$$k_*=k_*\Big(R, \|(\vm,\wm)\|_{L^\infty(0,+\infty)}, D_*, \|(y,D)\|_{C^1([0,+\infty)\times\SSS^2)} \Big)>0$$
such that, for every
$\pwh^0, \pvh^0 \in \mathbb{R}^3$ with $\|(\pwh^0-\pw, \pvh^0-\pv)\|\leqslant R$, $k_y, k_D>k_*$,
there exists a unique solution
$(\yh,\Dh,\pwh,\pvh) \in C^1([0,+\infty)\times\SSS^2,\mathbb{R})^2 \times C^1([0,+\infty),\mathbb{R}^3)^2$
of the Cauchy-problem
\begin{equation} \label{Observer_CYpb}
\left\lbrace \begin{array}{ll}
 \partial_t \yh &= - \nabla \yh \cdot(\etab \times (\wm-\pwh +\frac{1}{D}\etab\times (\vm-\pvh))) +k_y(y-\yh),
\\
 \partial_t \Dh &= - \nabla \Dh \cdot(\etab \times (\wm-\pwh +\frac{1}{D}\etab\times (\vm-\pvh))) -(\vm-\pvh)\cdot\etab \\
 &\quad+k_D(D-\Dh),
\\
\frac{d \pwh}{dt} &= -\int_{\SSS^2}\left( \lambda_y(\yh-y)(\nabla \yh\times\etab) +\lambda_D(\Dh-D)(\nabla \Dh\times\etab) \right)d\sigma_\eta,
\\
\frac{d \pvh}{dt} &= -\int_{\SSS^2}\left( \lambda_D(\Dh-D)\etab+\lambda_y(\yh-y)(\frac{1}{D}\etab\times(\etab\times\nabla \yh)\right.\\
&\quad\left.+\lambda_D(\Dh-D)(\frac{1}{D}\etab\times(\etab\times\nabla \Dh)  \right)d\sigma_\eta,
\\
\yh(0,\etab)&=y(0,\etab),
\\
\Dh(0,\etab)&=D(0,\etab),
\\
\pwh(0)&=\pwh^0,
\\
\pvh(0)&=\pvh^0.
\end{array} \right.
\end{equation}
\end{proposition}

The proof of Proposition \ref{Prop:WP1} is technical, thus we postpone it to Appendix \ref{App:WP}.

\subsection{Regularity and bounds}
The following proposition states the uniform boundedness of the estimates and their partial derivatives.

\begin{proposition} \label{Prop:WP2}
Let $R$ and $k_*$ be as in Proposition \ref{Prop:WP1}.
There exists
$$\mathcal{C}=\mathcal{C}\Big( R, \|(\vm,\wm)\|_{L^\infty(0,+\infty)}, D_*, \|(y,D)\|_{C^1([0,+\infty)\times\SSS^2)} \Big)>0$$
such that for every $(\pwh^0, \pvh^0) \in \mathbb{R}^6$ with $\|(\pwh^0-\pw, \pvh^0-\pv)\| \leqslant R$, $k_y, k_D>k_*$,
the solution of the Cauchy problem (\ref{Observer_CYpb}) satisfies
$$\left\| \yh(t,\etab)\right\|,\|\Dh(t,\etab)\|,\|\nabla \yh(t,\etab)\|, \|\nabla \Dh(t,\etab)\| \leqslant \mathcal{C}, \quad \forall (t,\etab) \in [0,+\infty)\times \SSS^2.$$
\end{proposition}
The proof of this result is a corollary of an intermediary result of the proof of Proposition \ref{Prop:WP1}
(see Proposition \ref{Prop:borneC1_pr_iteration} in Appendix).

\subsection{Continuity}
The following proposition details the continuity of the estimated state with respect to $(y,D,\om,\vit)$.

\begin{proposition}\label{Prop:WP3}
Let $T>0$.
Let $(\vit^n,\om^n)_{n\in \mathbb{N}}$, $(\vit,\om) \in C^1_b([0,+\infty),\mathbb{R}^3)$,
$(y^n,D^n)_{n \in \mathbb{N}}$, $(y,D) \in C^1_b([0,+\infty)\times \SSS^2,\mathbb{R})$ be associated solutions of (\ref{eq:state_equations})
such that
\begin{equation} \label{hyp:CV(v,w)}
(\vit^n,\om^n) \underset{ n \rightarrow + \infty}{\longrightarrow} (\vit,\om) \quad \text{ in } C^0([0,T],\mathbb{R}^3)^2,
\end{equation}
\begin{equation} \label{hyp:bound(y,G)}
(y^n,D^n)_{n \in \mathbb{N}} \quad \text{ is bounded in } C^1([0,T] \times \SSS^2,\mathbb{R})^2,
\end{equation}
\begin{equation} \label{hyp:CV(y,G)}
(y^n,D^n) \underset{ n \rightarrow + \infty}{\longrightarrow} (y,D) \quad \text{ in } C^0([0,T] \times \SSS^2,\mathbb{R})^2.
\end{equation}
Let $R$ and $k_*$ be as in Proposition \ref{Prop:WP1},
$(\yh^n,\Dh^n,\pwh^n,\pvh^n)_{n \in \mathbb{N}}$, $(\yh,\Dh,\pwh,\pvh)$ be associated solutions of the observer (\ref{Observer_CYpb})
with $k_y, k_D>k_*$ and initial conditions such that
$$(\pwh^n,\pvh^n)(0) \rightarrow (\pwh,\pvh)(0) \quad \text{ when } n \rightarrow + \infty.$$
Then, for every $T>0$,
\begin{equation} \label{CV(pwh,pvh)}
(\pwh^n,\pvh^n) \underset{ n \rightarrow + \infty}{\longrightarrow} (\pwh,\pvh) \quad \text{ in } C^0([0,T],\mathbb{R}^3)^2,
\end{equation}
\begin{equation} \label{CV(yh,Gh)}
(\yh^n,\Dh^n)(t) \underset{ n \rightarrow + \infty}{\longrightarrow} (\yh,\Dh)(t) \quad \text{ in } C^0([0,T] \times \SSS^2,\mathbb{R})^2, \quad \forall t \in [0,T].
\end{equation}
\end{proposition}
The proof of Proposition \ref{Prop:WP3} relies on tools introduced in the proof of Proposition \ref{Prop:WP1}, thus, we postpone it to Appendix \ref{sec:app_cont}.

\subsection{Convergence}\label{subsec:convergence}
The following theorem states that the proposed observer is optimal, in the sense that the conditions for convergence of the observer and observability of the system are strictly equivalent.

\begin{theorem} \label{thm:CV}
Let $R$ and $k_*$ be as in Proposition \ref{Prop:WP1}.
For every $(\pwh^0, \pvh^0) \in \mathbb{R}^6$ with $\|(\pwh^0-\pw, \pvh^0-\pv)\| \leqslant R$, $k_y, k_D>k_*$,
the solution of the Cauchy problem (\ref{Observer_CYpb}) satisfies
\begin{remunerate}
\item $\|(\yh-y)(t)\|_{L^\infty(\SSS^2)} + \|(\Dh-D)(t)\|_{L^\infty(\SSS^2)} \underset{t \rightarrow +\infty}{\longrightarrow} 0$\,,
\item for every adherence value $(\overline{\pv},\overline{\pw})$ of $(\pvh,\pwh)$ as $[t \rightarrow + \infty]$,
there exists $\bar C \in K \subset \Omega$  a camera position such that, $\forall \etab \in \SSS^2$
\begin{equation}
\nabla \bar{y} \cdot \etab \times((\overline{\pw}-\pw)+\frac{1}{\bar{D}} \etab \times (\overline{\pv}-\pv))=0
\label{eq:invarianty}
\end{equation}
\begin{equation}
\nabla \bar{D} \cdot \etab \times((\overline{\pw}-\pw)+\frac{1}{\bar{D}} \etab \times (\overline{\pv}-\pv))+\etab\cdot (\overline{\pv}-\pv)=0
\label{eq:invariantGam}
\end{equation}
\end{remunerate}
where $\bar{y}(t,\etab)$ and $\bar{D}(t,\etab)$ are the brightness and depth fields associated to the position $\bar{C}$ of the camera.
According to Theorem~\ref{thm:IS},
when the system is observable (i.e. when $(\Sigma, y_\Sigma)$ does not admit a rotation axis),
then $(\pwh,\pvh)(t) \underset{t \rightarrow +\infty}{\longrightarrow} (\pw,\pv)$.
\end{theorem}

{\em Proof of Theorem \ref{thm:CV}:}

\textbf{Step 1: Convergence of $f(t)$, defined by (\ref{eq:Lyapdecrease}), to zero when $t \rightarrow + \infty$.}
The function $V(t)$ is nonincreasing and nonnegative, thus it converges when $t \rightarrow + \infty$.
From (\ref{eq:Lyapdecrease}), we deduce that the nonnegative fonction $f$ satisfies
$$\int_0^t f(s) ds \leqslant V(0)-V(t)\,, \quad \forall t >0\,.$$
Thus $f \in L^1(0,+\infty)$. In order to conclude, it is sufficient to prove that $f$ is uniformly continuous on $[0,+\infty)$ (Barbalat Lemma).
We have
\begin{multline}
\dot{f}=2\int_{\SSS^2}\left( \lambda_y^*\ty\left(\nabla \yh \cdot \etab \times ( \tpw +\frac{1}{D}\etab\times \tpv )-\nabla\ty\cdot\frac{1}{D}\etab\times(\etab\times \vit)-k_y\ty\right)\right.\\
\left.+\lambda_D^*\tD\left(\nabla \Dh \cdot(\etab \times (\tpw +\frac{1}{D}\etab\times \tpv))-\etab\cdot\tpv-\nabla\tD\cdot\frac{1}{D}\etab\times(\etab\times \vit)-k_D\tD\right)
\right) d\sigma_\eta
\end{multline}
with $\lambda_y^*:=\lambda_y\left(k_y-\frac{L}{2}\right)$ and $\lambda_D^*:=\lambda_D\left(k_D-\frac{L}{2}\right)$.
Thanks to Proposition \ref{Prop:WP2}, we get $\dot{f} \in L^\infty(0,+\infty)$, which ends the proof of Step 1.
\\

\textbf{Step 2: Proof of statement 1.}
The function $(\yh-y,\Dh-D)$ belongs to $C^1_b([0,+\infty)\times \SSS^2)$ thanks to Proposition \ref{Prop:WP2} and (\ref{yDC1b}).
Thus Ascoli theorem guarantees the existence of adherence values in $C^0(\SSS^2,\mathbb{R})^2$ of $(\yh-y,\Dh-D)(t)$
when $t \rightarrow + \infty$. Thanks to Step 1, the only possible adherence value is $(0,0)$,
thus the whole function converges to $(0,0)$ in $C^0(\SSS^2,\mathbb{R})^2$.
\\

\textbf{Step 3: Proof of statement 2.}
Let $(\overline{\pv}^0,\overline{\pw}^0)$ be an adherence value of $(\pvh,\pwh)$ as $t \rightarrow + \infty$.
Let $(t_n)_{n \in \mathbb{N}}$ be an increasing sequence of $[0,+\infty)$ such that
$(\pvh,\pwh)(t_n) \rightarrow (\overline{\pv}^0,\overline{\pw}^0)$.
Let $T>0$.
The functions $\vit,\om$ (resp. $y, D$) belong to $C^1_b([0,+\infty),\mathbb{R}^3)$ (resp. $C^1_b([0,+\infty)\times\SSS^2,\mathbb{R})$),
thus Ascoli's theorem guarantees the existence of
$\overline{\vit}, \overline{\om} \in C^0([0,T],\mathbb{R}^3)$
(resp. $\overline{D}^0, \overline{y}^0 \in C^0(\SSS^2,\mathbb{R})$)
such that, up to an extraction
$$\begin{array}{c}
(\vit,\om)(t_n+.) \rightarrow (\overline{\vit},\overline{\om})(.) \text{ in } C^0([0,T],\mathbb{R}^3),\\
$$(y,D)(t_n) \rightarrow (\overline{y}^0, \overline{D}^0) \text{ in } C^0(\SSS^2,\mathbb{R}).
\end{array}$$
The continuity of the flow for \eqref{eq:state_equations} justifies that
$$(y,D)(t_n+.) \rightarrow (\overline{y},\overline{D})(.) \text{ in } C^0([0,T]\times\SSS^2,\mathbb{R})^2$$
where
$$\left\lbrace \begin{array}{l}
\partial_t \overline{y} = -  \nabla \overline{y} \cdot(\etab \times (\overline{\om} +\frac{1}{\overline{D}}\etab\times \overline{\vit})), \\
\partial_t \overline{D} = - \nabla \overline{D} \cdot(\etab \times (\overline{\om} +\overline{D}\etab\times \overline{\vit})) - \overline{\vit}\cdot\etab,\\
\overline{y}(0)=\overline{y}^0,\\
\overline{D}(0)=\overline{D}^0.
\end{array} \right.$$
Thanks to Proposition \ref{Prop:WP3}, we know that
$$\begin{array}{c}
(\pwh,\pvh)(t_n+.) \rightarrow (\widehat{\overline{\pw}},\widehat{\overline{\pv}}) \quad \text{ in } C^0([0,T],\mathbb{R}^3)^2,\\
$$(\yh,\Dh)(t_n+.) \rightarrow (\widehat{\overline{y}},\widehat{\overline{D}}) \quad \text{ in } C^0([0,T]\times\SSS^2,\mathbb{R})^2,
\end{array}$$
where $(\widehat{\overline{y}},\widehat{\overline{D}}, \widehat{\overline{\pw}},\widehat{\overline{\pv}})$
is the solution of the observer associated to $(\overline{y},\overline{D})$,
$\overline{\om}_m(t)=\overline{\om}(t)+\pw$, $\overline{\vit}_m(t)=\overline{\vit}(t)+\pv$
and the initial conditions
$(\widehat{\overline{\pw}},\widehat{\overline{\pv}})(0)=(\overline{\pw}^0,\overline{\pv}^0)$.

We know that $V(t)$ converges to some limit $V_\infty$ as $t \rightarrow + \infty$.
Moreover, $V(t_n+t) \rightarrow \overline{V}(t)$ for every $t \in [0,T]$.
Therefore, $\overline{V}(t)=V_\infty$ for every $t \in [0,T]$.
In particular, $0 \leqslant \overline{f}(t) \leqslant |d\overline{V}/dt| \equiv 0$ i.e.
$(\widehat{\overline{y}}-\overline{y})(t)=(\widehat{\overline{D}}-\overline{D})(t)=0, \forall t \in [0,T]$.
Substracting the equations on $\widehat{\overline{y}}$ and $\overline{y}$ (resp. on $\widehat{\overline{D}}$
and $\overline{D}$), we get the relation (\ref{eq:invarianty}) (resp. (\ref{eq:invariantGam})). \hfill $\Box$

\section{Practical implementation and simulations}
\label{sec:implem}
\subsection{Adaptation to a spherical cap}
Concretely, a spherical camera is only a model, and the image perceived by real cameras only occupy a part of $\SSS^2$. Let us call $K$ this portion: $y(t,\etab)$ and $D(t,\etab)$ are known only for $\etab \in K$. The  observer introduced in \ref{sec:asymp_obs} can not be readily used  since it brings into play the integral of $y$ or $D$ over the whole unit sphere. We will see that one can compensate this problem by considering virtual observations, equal to the real observations over the window defined by $K$.
Let $K_1$ and $K_2$ be two compact sets s.t. $\stackrel{\circ}{K_1} \subset \stackrel{\circ}{K_2} \subset \stackrel{\circ}{K} $. Let $\phi$ be a ${\cal C}^{\infty}$ scalar field $\mathbb S^2\ni \etab \mapsto \phi(\etab)\in\R$, s.t.
$\phi=1$ on  $K_1$, $\phi=0$ and $\nabla\phi=0$ on $\SSS^2 \smallsetminus \stackrel{\circ}{K_2}$.
Let us define $X=\phi y$ and $\Lambda=\phi D$.
Then,
$$
\partial_t X =-\nabla X \cdot(\etab \times (\om +\frac{1}{D}\etab\times \vit))+y\nabla \phi \cdot(\etab \times (\om +\frac{1}{D}\etab\times \vit))
$$
and
$$
\partial_t\Lambda =-\nabla \Lambda \cdot(\etab \times (\om +\frac{1}{D}\etab\times \vit))+D\nabla \phi \cdot(\etab \times (\om +\frac{1}{D}\etab\times \vit))-\phi\etab\cdot \vit
.
$$
We propose the following adaptation of observer~\eqref{eq:observer}:
\begin{equation}
\left\lbrace
    \begin{array}{ll}
        \partial_t\Xh=&- \nabla \Xh \cdot(\etab \times (\wm-\pwh +\frac{1}{D}\etab\times (\vm-\pvh)))
        \\
        &+y\nabla \phi \cdot(\etab \times(\wm-\pwh +\frac{1}{D}\etab\times (\vm-\pvh)))
        \\
        & +k_X(X-\Xh),\\
        \partial_t\Lh=&- \nabla \Lh \cdot(\etab \times (\wm-\pwh +\frac{1}{D}\etab\times (\vm-\pvh)))\\
        &+D\nabla \phi \cdot(\etab \times(\wm-\pwh +\frac{1}{D}\etab\times (\vm-\pvh)))\\
        &-\phi(\vm-\pvh)\cdot\etab +k_{\Lambda}(\Lambda-\Lh),\\
        \partial_t\pwh=&-k_{\omega}\int_{K}\lambda_X(\Xh-X)((\nabla \Xh-y\nabla\phi)\times\etab)\\
        &+ \lambda_{\Lambda}(\Lh-\Lambda)((\nabla \Lh-D\nabla\phi)\times\etab) d\sigma_\eta,\\
        \partial_t\pvh=&-k_{v}\int_{K} \lambda_{\Lambda}(\Lh-\Lambda)\etab\\
        &+\lambda_X(\Xh-X)(\frac{1}{D}\etab\times(\etab\times(\nabla \Xh-y\nabla\phi))\\
        &+\lambda_{\Lambda}(\Lh-\Lambda)(\frac{1}{D}\etab\times(\etab\times(\nabla \Lh-D\nabla\phi)) d\sigma_\eta.
    \end{array}
    \right.
\label{eq:obs_reduit}
\end{equation}
Notice that the estimates $\Xh$ and $\Lh$ of the virtual observations are computed only on the domain $K$, and not on the entire sphere $\SSS^2$, and only their values on $K$ are used in the dynamics of the estimated biases. Neumann boundary conditions are imposed on the parts of $\partial K$ where the direction of propagation points toward the inside of the domain $K$:
$$\frac{\partial \Xh}{\partial \mathbf{n}}=0,\quad \frac{\partial \Lh}{\partial \mathbf{n}}=0$$
where $\mathbf n$ is the outwards-pointing normal of $\partial K$ such that $$\mathbf{n}\cdot (\etab \times (\wm-\pwh +\frac{1}{D}\etab\times (\vm-\pvh)))<0.$$
Let us choose the candidate Lyapunov function
$$
V=\frac{1}{2}\left(\int_{K}\left(\lambda_X\tX^2+\lambda_{\Lambda}\tL^2 \right)d\sigma_\eta + \frac{\tpw^2}{k_{\omega}} + \frac{\tpv^2}{k_v}  \right).
$$
One can prove that

\begin{multline*}
\dot{V}=\int_{K} \left(\lambda_X \tX\left(- \nabla \tX \cdot(\etab \times (\om +\frac{1}{D}(\etab\times \vit)))-k_X\tX\right)\right)d\sigma_\eta\\
+\int_{K} \left(\lambda_{\Lambda}\tL\left(- \nabla \tL \cdot(\etab \times (\om +\frac{1}{D}(\etab\times \vit)))-k_{\Lambda}\tL\right)\right)d\sigma_\eta
\end{multline*}
and integration by parts yields
\begin{multline*}
\dot{V}=-\int_{K} \left(\lambda_X\left(k_X-\frac{\nabla\cdot W(t,\etab)}{2}\right)\tX^2+ \lambda_{\Lambda}\left(k_{\Lambda}-\frac{\nabla\cdot W(t,\etab)}{2}\right)\tL^2\right)d\sigma_\eta\\
-\oint_{\partial K}\left(\lambda_X\tX^2+\lambda_{\Lambda}\tL^2\right)\frac{ W(t,\etab)}{2}\cdot\bold n dl_\eta
\label{eq:lyapunov_calotte}
\end{multline*}
where $W$ is defined as $(\R,\SSS^2)\ni(t,\etab)\mapsto W(t,\etab)=\etab\times(\omega+\frac{1}{D}\etab\times v)\in\R^3$.
We guess that for given environment and trajectory, sufficiently large correction gains $k_X$ and $k_{\Lambda}$ can ensure the decrease of the candidate Lyapunov function $V$. Then convergence analysis done when $K=\SSS^2$ should be extended  to compact sub-domains $K$ of  $\SSS^2$. As the candidate Lyapunov function should decrease unless $\tX=\tL=0$ on $K$, we guess that the necessary and sufficient condition for observability is the same as stated in Theorem \ref{thm:IS}, but restricted to the visible part of the environment.

\subsection{The observer in pinhole coordinates}
The previous observer~\eqref{eq:obs_reduit} can be finally adapted to a real model of camera: we choose here the widely spread pinhole camera model enabling a correspondence between the local coordinates on $\SSS^2$ with a rectangular grid of pixels. The pixel of coordinates $(z_1,z_2)$ corresponds to the unit vector $\etab\in\SSS^2$ of coordinates in $\R^3$:  $\left(1+z_{1}^2+z_{2}^2\right)^{-1/2}(z_1,z_2,1)^T$. The optical camera axis (pixel $(z_1,z_2)=(0,0)$) corresponds here  to the direction $z_3$. Directions $1$ and $2$ correspond respectively to the horizontal axis from left to right and to the vertical axis from top to bottom on the image frame.

The gradients $\nabla y$ and $\nabla D$ must be expressed with respect to $z_1$ and $z_2$. Let us detail this derivation for $y$. Firstly, $\nabla y$ is tangent to $\mathbb S^2$, thus $\nabla y\cdot\etab=0$. Secondly, the differential $dy$ corresponds to $\nabla y\cdot d\etab$ and to $\dv{y}{z_{1}} dz_{1} +\dv{y}{z_{2}} dz_{2}$. By identification, we get the Cartesian coordinates of $\nabla y$ in $\R^3$. Similarly we get the three coordinates of $\nabla D$. Plugging these expressions in~\eqref{eq:obs_reduit}, we get a partial differential equations (PDE) system written in terms of $(t,z_1,z_2)$ as independent variables. Due to space limitation, this system is not given here, but its derivation is straightforward and a little bit tedious.

\subsection{Simulations}
\begin{figure}
      \centering
      \includegraphics[scale=0.53]{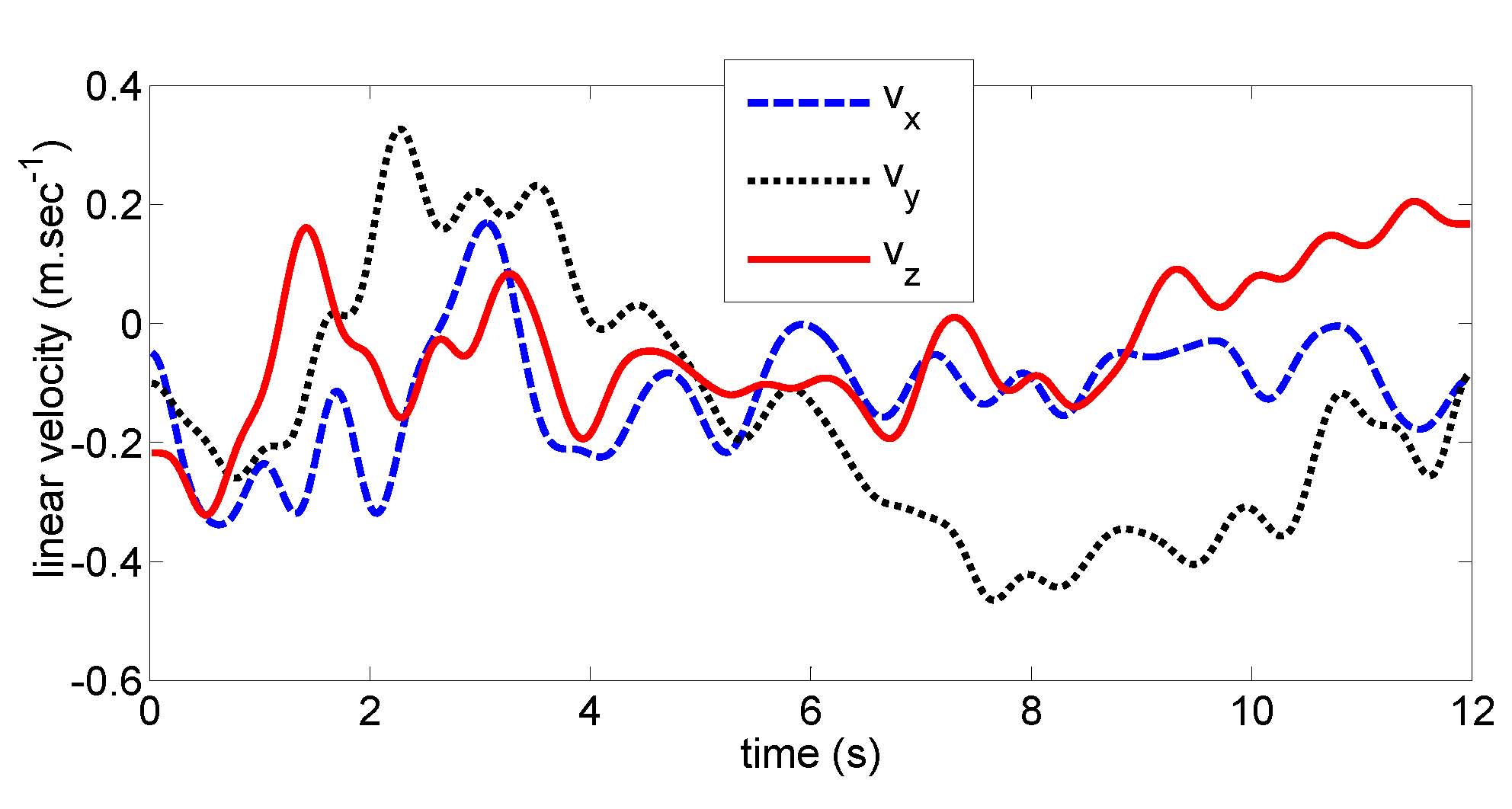}
      \caption{Components of the linear velocity $v$ used to test the observer: translations in the horizontal, vertical and optical axis directions, respectively.}
      \label{fig:linvelocity}
      \includegraphics[scale=0.53]{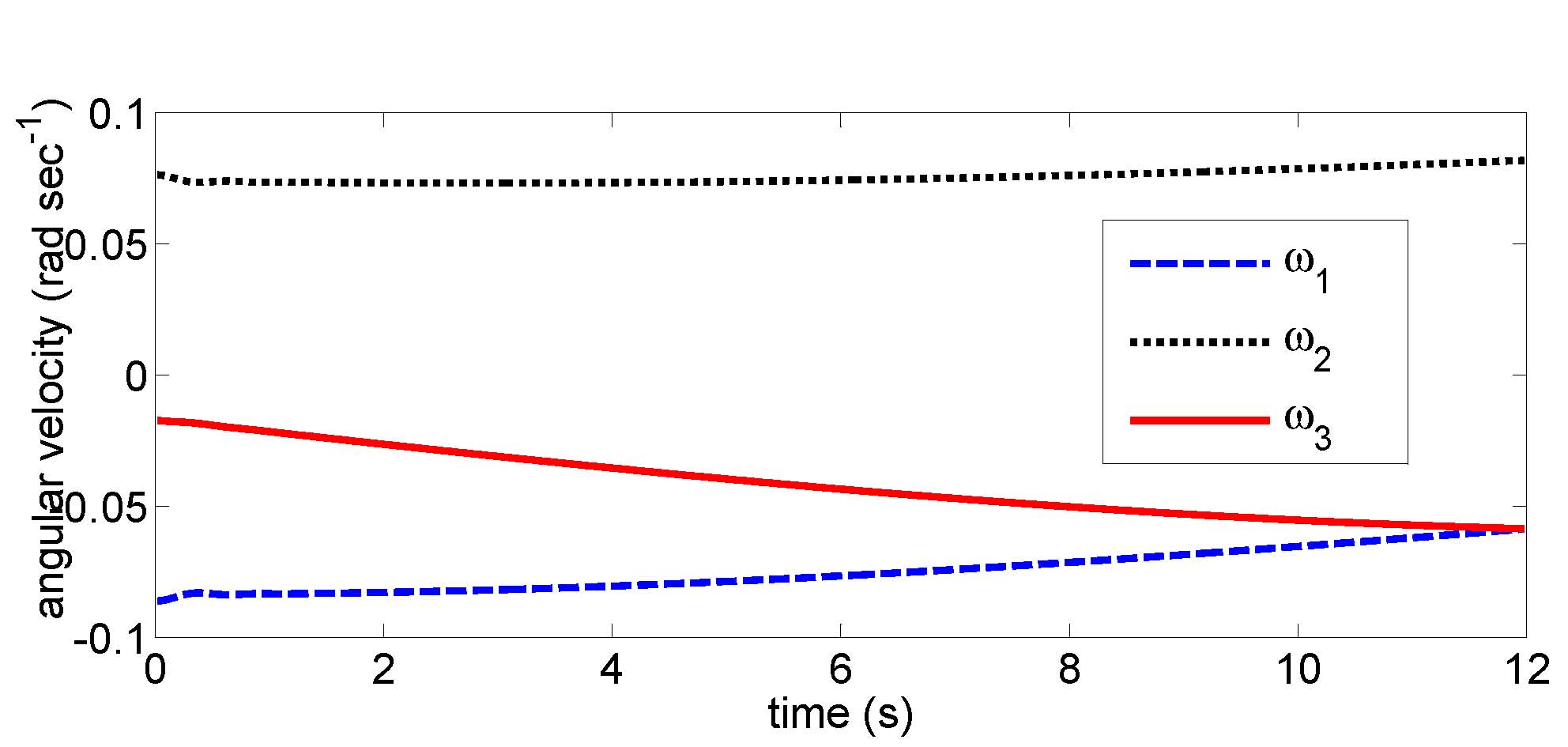}
      \caption{Components of the angular velocity $\omega$ used to test the observer: rotations around the horizontal, vertical and optical axis, respectively.}
      \label{fig:angvelocity}
\end{figure}

The non-linear asymptotic observer~\eqref{eq:obs_reduit}  is tested on a sequence of synthetic images characterized by the following:
\begin{romannum}
\item \textit{virtual camera providing images  restricted to  the spherical cap $K\subset\SSS^2$}: the size of each image is 640 by 480 pixels, the frame rate of the sequence is 42 Hz and the field of view is 50 deg by 40 deg.
\item \textit{motion of the virtual camera $(\om,\vit)$}: it consists of the motion of a real hand-held camera (filtered data), combining translations and rotations in each direction; the real linear and angular velocities expressed in the camera frame are plotted in Fig.\ref{fig:linvelocity} and Fig.\ref{fig:angvelocity}; zero-mean normally distributed noise (standard deviations $\sigma_v$ and $\sigma_{\omega}$) is added to these velocities to test the robustness;
\item \textit{virtual scene corresponding to $(\Sigma,y_\Sigma)$}: it consists of the walls, ceiling and floor of a virtual room; the observed   walls  are virtually painted with a gray pattern, whose intensity varies in horizontal and vertical directions as a sinusoid function;
\item \textit{generation of the brightness images $y$}: each pixel of an image has an integer value varying from $1$ to $256$, directly depending on the intensity of the observed surface in the direction indexed by the pixel, to which a zero-mean normally distributed noise with standard deviation $\sigma_y$ is added to test the robustness;
\item \textit{generation of the depth images $D$}: to each pixel of the rectangular grid of an image is attributed the depth of the corresponding element of the observed surface, computed with respect to position and orientation of the camera in the room, to which a zero-mean normally distributed noise with standard deviation $\sigma_D$ is added to test the robustness.
\end{romannum}

\begin{figure}
      \centering
     \includegraphics[scale=0.53]{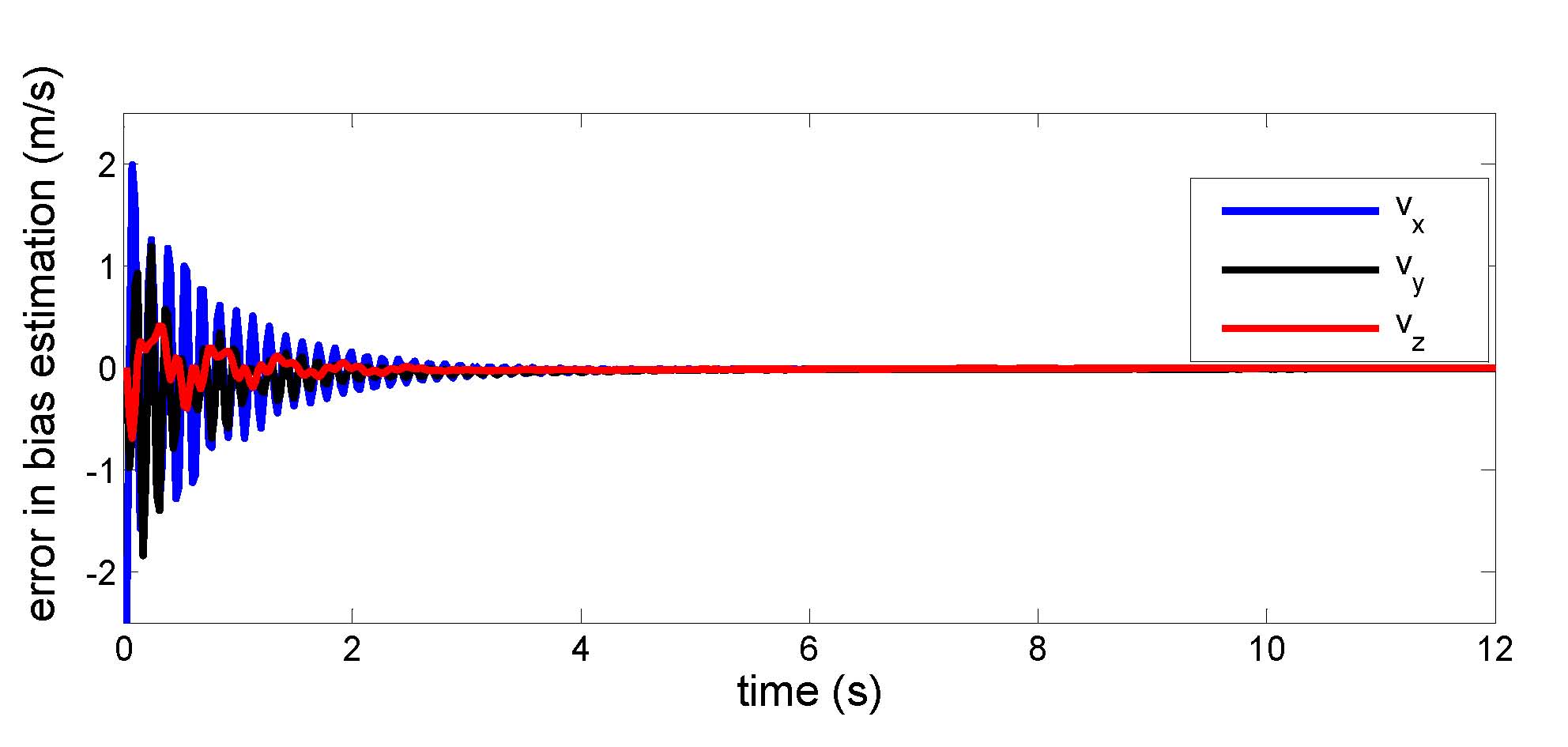}
      \caption{Error of estimation of the biases in translation for noiseless  image and depth data. Real biases are $p_{\omega_1}=0.05$ rad.$s^{-1}$ and $p_{v_x}=2.5$ m.$s^{-1}$}
      \label{fig:biaistrans}
     \includegraphics[scale=0.53]{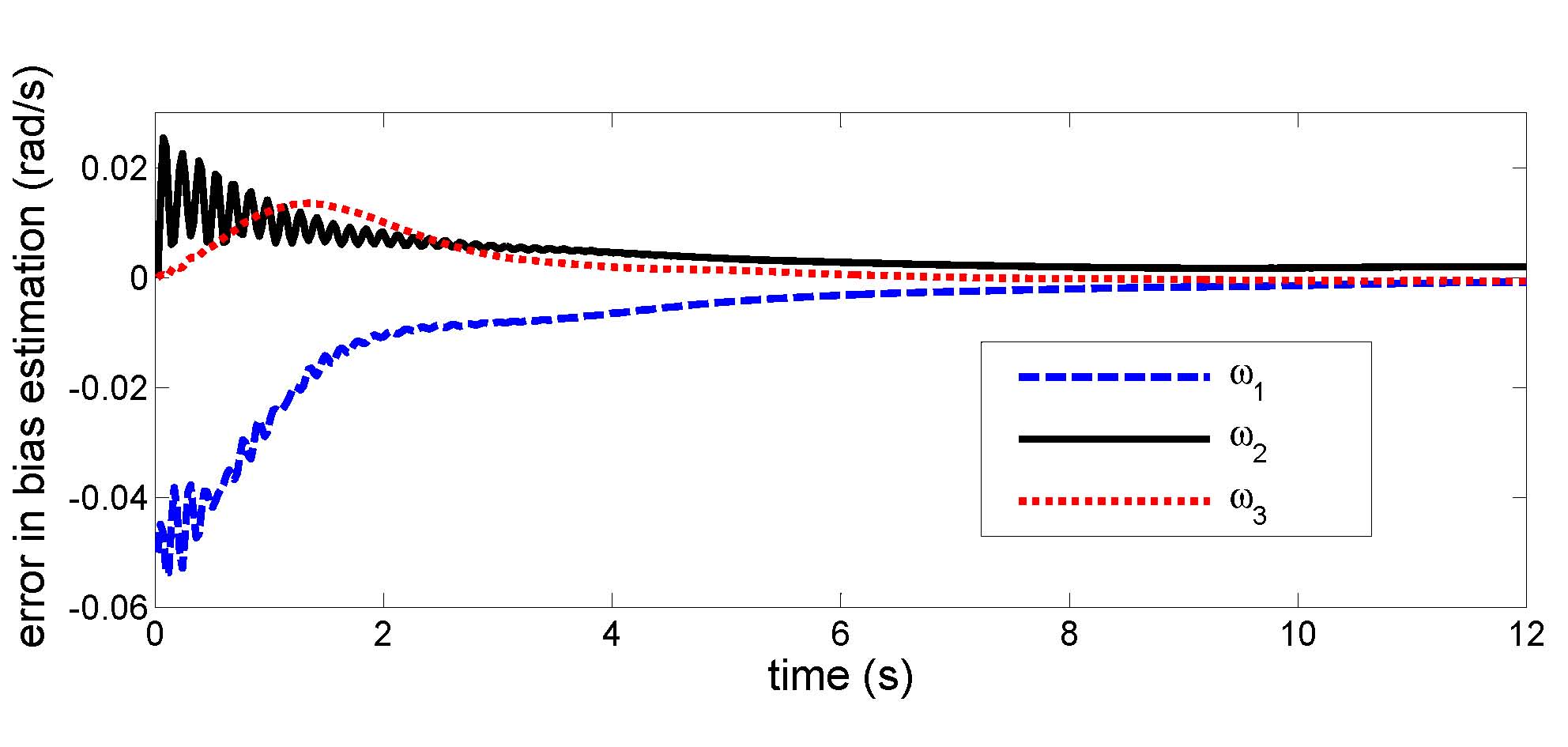}
      \caption{Error of estimation of the biases in rotation for noiseless image and depth data. Real biases are $p_{\omega_1}=0.05$ rad.$s^{-1}$ and $p_{v_x}=2.5$ m.$s^{-1}$}
      \label{fig:biaisrot}
\end{figure}

The numerical resolution used to compute $\yh$, $\Dh$, $\pwh$ and $\pvh$ according to \eqref{eq:obs_reduit} is based on a temporal Euler discretization scheme where $\nabla \yh$ and $\nabla \Dh$ are computed via differentiation filters (Sobel filtering) directly from the image and depth previous estimates.
The observer is then tested for reasonable biases: in rotation, a bias of $0.05$ rad.$s^{-1}$ ($10,000$ deg/h, for a low-cost gyroscope) around the horizontal axis; in translation, a bias of $2.5$ m.$s^{-1}$ in the horizontal direction ($9$ km/h, for the windspeed). In other words, $p_{\omega_1}=0.05$ rad.$s^{-1}$ and $p_{v_x}=2.5$~m.$s^{-1}$. Biases in the other directions are set to 0. Initial conditions for $\yh$ and $\Dh$ are $y(0,\etab)$ and $D(0,\etab)$. Initial conditions for the estimated biases are set to zero. The chosen correction gains are: $k_y=k_{D}=2 s^{-1}$, $k_v=10^{-2} \text{m}^2.s^{-2}$ and $k_{\omega}=10^{-5} \text{rad}.\text{m}.s^{-2}$.  These correction gains are chosen in accordance with the expected values of biases and the scene averaged depth, to enable a reasonable convergence speed. The correction gains $k_y$ and $k_{D}$ are comparatively much larger than $k_v$, which is itself larger than $k_{\omega}$, as large oscillations in the estimation of $p_{\omega}$ can make the discretized observer to diverge. Finally, the ponderation coefficients are $\lambda_y=1$ and $\lambda_{D}=5000 \text{m}^2$, chosen to compensate the difference of magnitudes of $y$ and $D$. First, when image and depth data contain no noise ($\sigma_y=\sigma_D=0$), the results are plotted in Fig.\ref{fig:biaistrans} and Fig.\ref{fig:biaisrot} as the instantaneous errors of estimations $\tpv$ and $\tpw$ expressed in the camera frame, respectively. In the first $6$ s, errors slowly converge towards $0$, and coupling between rotation and translation occurs: this reflects the fact that an horizontal translation can be interpreted as a rotation around the vertical axis to a certain extent. Oscillations decrease, and after convergence, errors stay bounded: for the bias in rotation, it does not exceed $0.002$ rad.$s^{-1}$ ($4$ \% of the original bias); in translation, it is less than $0.01$ m.$s^{-1}$ ($0.4$ \% of the original bias).

\begin{figure}
      \includegraphics[scale=0.53]{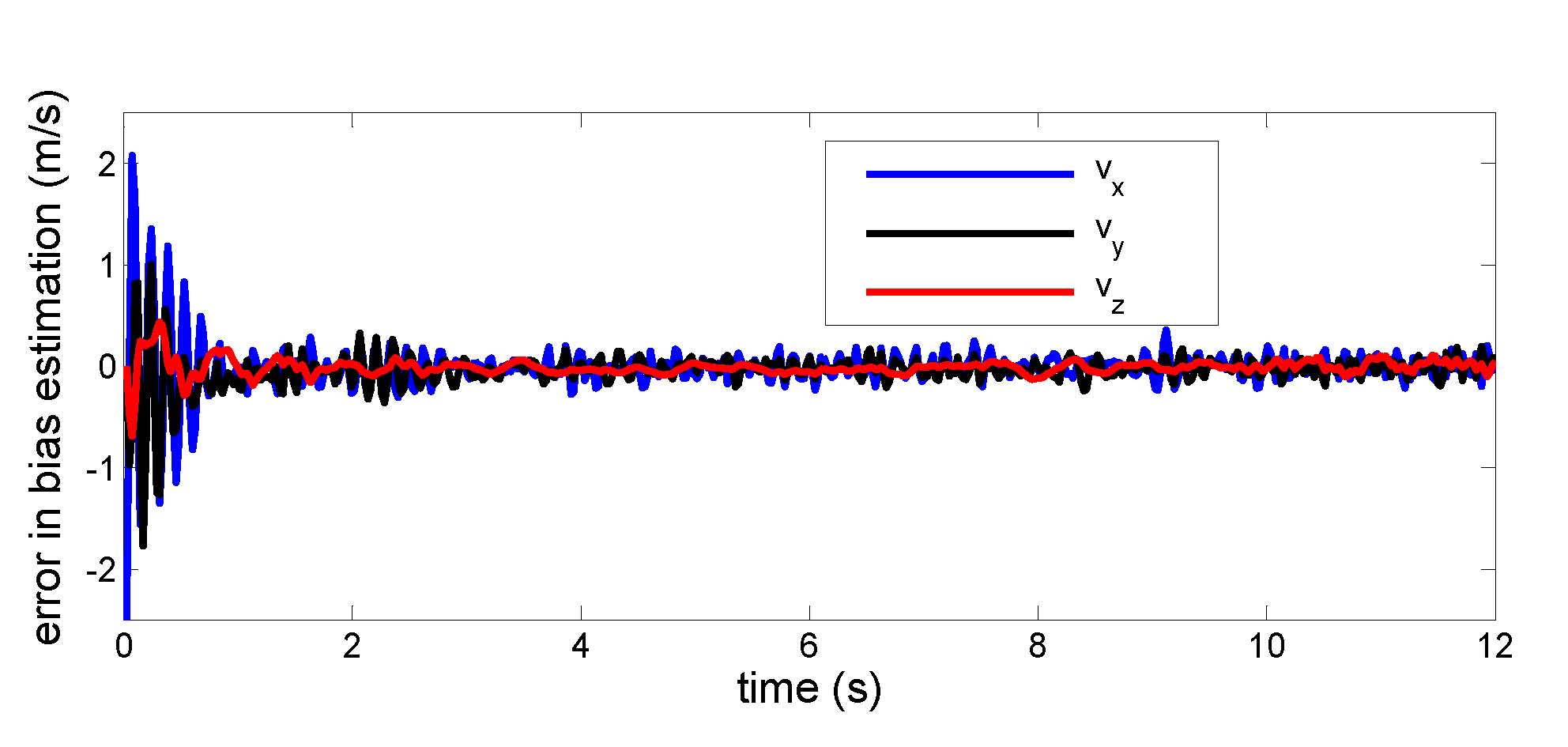}
      \caption{Error of estimation of the biases in translation for noisy image and depth data, and noisy velocities: $\sigma_y=30$, $\sigma_D=25$ cm, $\sigma_v=0.05$m.$s^{-1}$, $\sigma_{\omega}=0.005$ rad.$s^{-1}$. Real biases are $p_{\omega_1}=0.05$ rad.$s^{-1}$ and $p_{v_x}=2.5$ m.$s^{-1}$}
      \label{fig:biaistransbruit}
      \includegraphics[scale=0.53]{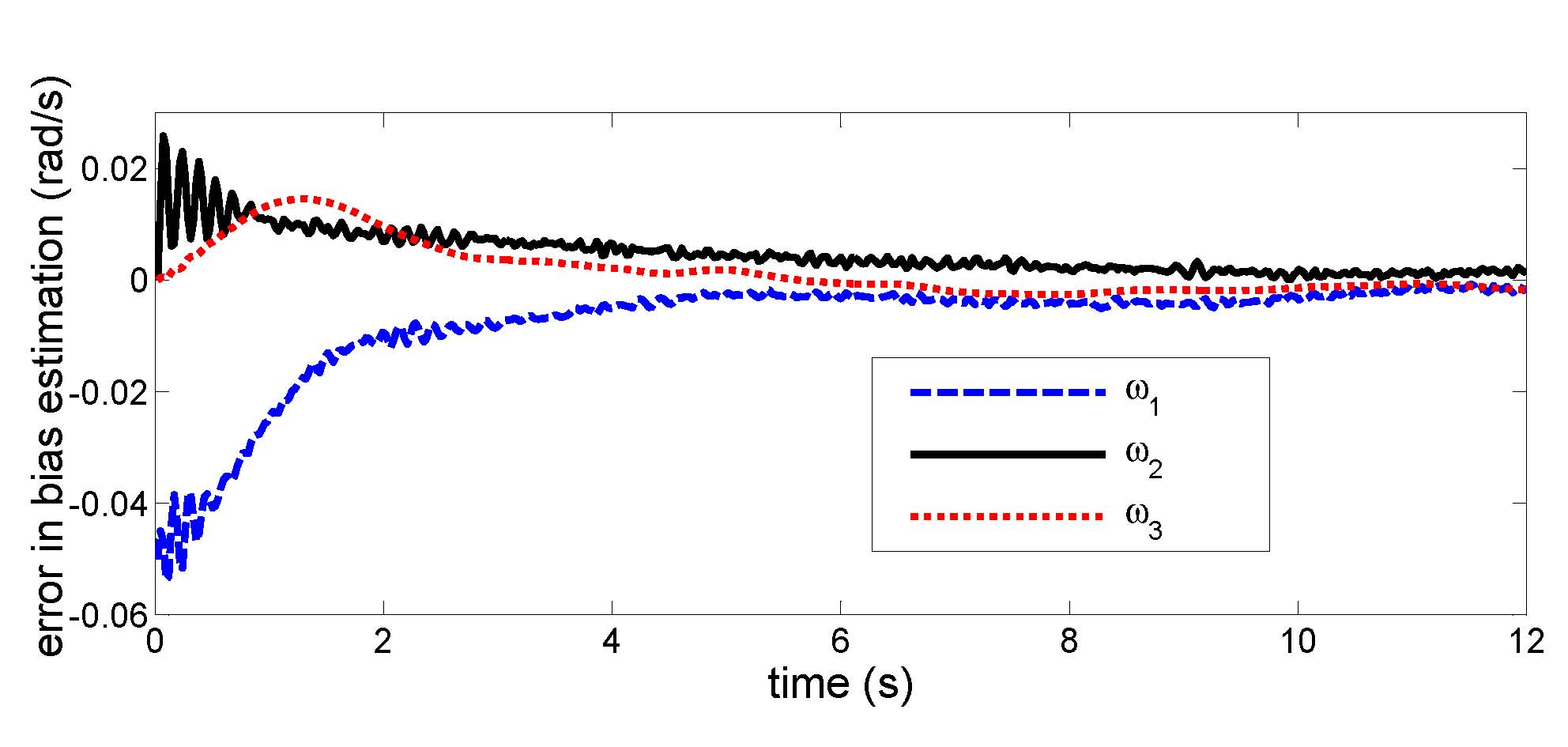}
      \caption{Error of estimation of the biases in rotation for noisy image and depth data, and noisy velocities: $\sigma_y=30$, $\sigma_D=25$ cm, $\sigma_v=0.05$m.$s^{-1}$, $\sigma_{\omega}=0.005$ rad.$s^{-1}$. Real biases are $p_{\omega_1}=0.05$ rad.$s^{-1}$ and $p_{v_x}=2.5$ m.$s^{-1}$}
      \label{fig:biaisrotbruit}
\end{figure}

Then, to test the robustness of the method, noise is added to the image data ($\sigma_y=30$, about $12\%$ of the full scale), to the depth data ($\sigma_D=25$ cm, which is three times as much that can be expected from a Kinect device), to linear velocity ($\sigma_v=0.05$m.$s^{-1}$) and to angular velocity ($\sigma_{\omega}=0.005$ rad.$s^{-1}$). Results are plotted in Fig.\ref{fig:biaistransbruit} and Fig.\ref{fig:biaisrotbruit}. Convergence time is shorter for biases in translation estimation: after 3 s, error does not exceed 0.2 m.$s^{-1}$ ($8$ \% of the original bias). For the rotation, convergence is slower (as $k_{\omega}$ is smaller), but in the last 3 sec of the simulation, biases are estimated up to $0.003$ rad.$s^{-1}$ ($6$ \% of the original bias).

\clearpage

\section{Conclusion}

We have proposed new  infinite dimensional nonlinear observers \eqref{eq:observer} and~\eqref{eq:obs_reduit} that simultaneously filters image  and depth data and estimates constant biases on angular and translational velocities. Observer design is based  on a Lyapunov functional and convergence analysis has been done under $C^1$ regularity assumptions. With Theorems~\ref{thm:IS} and~\ref{thm:CV},   we  have proved asymptotic  convergence of observer~\eqref{eq:observer} if the  scene  $(\Sigma,y_\Sigma)$ does not admit any cylindrical symmetry.   Preliminary simulations indicate that our  convergence analysis could certainly be extended  to  more general situations  with weaker regularities.  It will be interesting  to test the nonlinear observer~\eqref{eq:obs_reduit} on real data, since  Kinect-like devices can provide the necessary image and depth dense fields.

\section{Appendix: Proof of the well posedness of the observer}
\label{App:WP}

The goal of this section is the proof of Proposition \ref{Prop:WP1} in 3 steps:
\begin{itemize}
\item the first step consists in proving the existence of local (in time) solutions,
\item the second step consists in proving that solutions stays bounded in a appropriate way, so that
\item we can iterate the local argument to get global (in time) solutions, in a third step.
\end{itemize}

\subsection{Local solutions}
The goal of this section is the proof of the following result.

\begin{proposition} \label{Prop:WP1_local}
Let $R$ and $M>0$. There exists
$$T=T (M,R,\|(\vm,\wm)\|_{L^\infty(0,+\infty)},D_*,\|(y,D)\|_{C^1([0,+\infty)\times\SSS^2)},k_y,k_D,\lambda_y,\lambda_D )>0$$
such that, for every $\pwh^0$, $\pvh^0 \in \mathbb{R}^3$, $\yh^0$, $\Dh^0 \in C^1(\SSS^2)$ with
\begin{equation} \label{def:R}
\| (\pwh^0-\pw, \pvh^0-\pv)\| \leqslant R\,,
\end{equation}
\begin{equation} \label{borne_C1_yD0}
\|\yh^0\|_{C^1(\SSS^2)}\,, \|\Dh^0\|_{C^1(\SSS^2)} \leqslant M\,,
\end{equation}
there exists a unique solution
$(\yh,\Dh,\pwh,\pvh) \in C^1([0,T]\times\SSS^2,\mathbb{R})^2 \times C^1([0,T],\mathbb{R}^3)^2$,
defined on the time interval $[0,T]$,
of the Cauchy-problem
\begin{equation} \label{Observer_CYpb_loc}
\left\lbrace \begin{array}{ll}
 \partial_t \yh &= - \nabla \yh \cdot(\etab \times (\wm-\pwh +\frac{1}{D}\etab\times (\vm-\pvh))) +k_y(y-\yh)\,,
\\
 \partial_t \Dh &= - \nabla \Dh \cdot(\etab \times (\wm-\pwh +\frac{1}{D}\etab\times (\vm-\pvh))) -(\vm-\pvh)\cdot\etab \\
 &\quad+k_D(D-\Dh)\,,
\\
\frac{d \pwh}{dt} &= -\int_{\SSS^2}\left( \lambda_y(\yh-y)(\nabla \yh\times\etab) +\lambda_D(\Dh-D)(\nabla \Dh\times\etab) \right)d\sigma_\eta,\,
\\
\frac{d \pvh}{dt} &= -\int_{\SSS^2}\left( \lambda_D(\Dh-D)\etab+\lambda_y(\yh-y)(\frac{1}{D}\etab\times(\etab\times\nabla \yh)\right.\\
&\quad\left.+\lambda_D(\Dh-D)(\frac{1}{D}\etab\times(\etab\times\nabla \Dh)  \right)d\sigma_\eta\,,
\\
\yh(0,\etab)&=\yh^0(\etab),
\\
\Dh(0,\etab)&=\Dh^0(\etab),
\\
\pwh(0)&=\pwh^0,
\\
\pvh(0)&=\pvh^0.
\end{array} \right.
\end{equation}

\end{proposition}

{\em Proof of Proposition \ref{Prop:WP1_local}:}
We prove the existence and uniqueness of the non-linear and non-local system of PDEs \eqref{Observer_CYpb_loc} by a standard fixed-point approach.
Our strategy consists in explicitly solving the two first PDEs of \eqref{Observer_CYpb_loc} with fixed functions $t \mapsto (\Pwh,\Pvh)(t)$, instead of $t \mapsto (\pwh,\pvh)(t)$.
Then, the explicit solutions $(\Yh(t),\DDh(t))$ are plugged in the third and fourth ODEs of \eqref{Observer_CYpb_loc}.
We conclude by applying the Banach fixed point theorem (see \cite{brezis}) on the map $(\Pwh,\Pvh) \mapsto (\pwh,\pvh)$.
\\

Let $R, M>0$, $\pwh^0$, $\pvh^0 \in \mathbb{R}^3$, and $\yh^0$, $\Dh^0 \in C^1(\SSS^2)$ such that (\ref{def:R}) and (\ref{borne_C1_yD0}) hold.
For any $T>0$, the functional space
$$E_T:=\{ (\Pwh,\Pvh) \in C^0([0,T],\mathbb{R}^6) ; \|(\Pwh,\Pvh)(t)- (\pw,\pv) \| \leqslant 2R, \forall t \in [0,T] \}$$
equipped with the $L^\infty((0,T),\mathbb{R}^6)$-norm is a Banach space.
To any element $(\Pwh,\Pvh) \in E_T$, we associate
the solutions $(\Yh, \DDh,\pwh,\pvh) \in C^1([0,T] \times \SSS^2)^2 \times C^1([0,T],\mathbb{R}^3)^2$ of the Cauchy problem
\begin{equation} \label{syst-fixedpoint}
\left\lbrace \begin{array}{ll}
 \partial_t \Yh &= - \nabla \Yh \cdot(\etab \times (\wm-\Pwh +\frac{1}{D}\etab\times (\vm-\Pvh))) +k_y(y-\Yh)\,,
\\
 \partial_t \DDh &= - \nabla \DDh \cdot(\etab \times (\wm-\Pwh +\frac{1}{D}\etab\times (\vm-\Pvh))) -(\vm-\Pvh)\cdot\etab +k_D(D-\DDh)\,,
\\
\frac{d \pwh}{dt} &= -\int_{\SSS}\left( \lambda_y(\Yh-y)(\nabla \Yh\times\etab) +\lambda_D(\DDh-D)(\nabla \DDh\times\etab) \right)d\sigma_\eta\,,
\\
\frac{d \pvh}{dt} &= -\int_{\SSS}\left( \lambda_D(\DDh-D)\etab+\lambda_y(\Yh-y)(\frac{1}{D}\etab\times(\etab\times\nabla \Yh)\right.\\
&\qquad\quad\left.+\lambda_D(\DDh-D)(\frac{1}{D}\etab\times(\etab\times\nabla \DDh)  \right)d\sigma_\eta\,,
\\
\Yh(0,\etab)&=\yh^0(\etab)\,,
\\
\DDh(0,\etab)&=\Dh^0(\etab)\,,
\\
\pwh(0)&=\pwh^0\,,
\\
\pvh(0)&=\pvh^0\,.
\end{array} \right.
\end{equation}
Let us emphasize that this Cauchy problem is well posed because $\Yh$ and $\DDh$ have explicit expressions along the characteristics
(see \cite{leveque1992numerical,serre1999systems} for a description of the method of characteristics).
\\

Indeed, let us introduce the flow $\Phi_t$ on the sphere $\SSS^2$ associated to the ODE
\begin{equation} \label{EDO}
\left\lbrace \begin{array}{l}
\frac{d\etab}{dt}=\etab \times \Big( \wm(t)-\Pwh(t) +\frac{1}{D(t,\etab)}\etab\times [\vm(t)-\Pvh(t)] \Big)\,,\\
\etab(0)=\etab_0\,.
\end{array}\right.
\end{equation}
This means that the solution of (\ref{EDO}) is $\etab(t)=\Phi_t(\etab_0)$, $\forall t,\etab_0$.
Note that the right hand side of this equation is continuous in $t$ and $C^1$ in $\etab$, thus the flow is uniquely defined.
This flow is defined for every $t \in [0,T]$ because no explosion is possible in finite time: it lives on the sphere $\SSS^2$.
Moreover, for every $t \in [0,T]$, $\etab \mapsto \Phi_t(\etab)$ is a $C^1$-diffeomorphism of the sphere $\SSS^2$, thus, there exists
$$C=C(R,\|(\wm,\vm)\|_{L^\infty(0,T)},D_*,\|D\|_{C^1([0,T]\times\SSS^2)})>0$$
such that
\begin{equation} \label{bound:diff_flow}
\left\| \frac{\partial \Phi_t}{\partial \etab}(\etab) \right\| , \left\| \frac{\partial \Phi_t^{-1}}{\partial \etab}(\etab) \right\| \leqslant C,
\quad \forall (t,\etab) \in [0,T] \times \SSS^2.
\end{equation}

With the characteristics method, we get the following explicit expressions
$$\Yh[t,\Phi_t(\etab)]=\yh^0(\etab)e^{-k_y t} + \int_0^t k_y y[\tau,\Phi_\tau(\etab)] e^{k_y(\tau-t)} d\tau,$$
$$\begin{array}{ll}
\DDh[t,\Phi_t(\etab)]=
&
\Dh^0(\etab)e^{-k_D t}
\\ &
+ \int_0^t \left(
k_D D(\tau,\Phi_\tau(\etab)) - (v_m-\Pvh)(\tau).\Phi_\tau(\etab)
\right) e^{k_D(\tau-t)}  d\tau,
\end{array}$$
or, equivalently
\begin{equation} \label{Yh:explicit}
\Yh[t,\etab]=\yh^0[0,\Phi_t^{-1}(\etab)]e^{-k_y t} + \int_0^t k_y y[\tau,\Phi_\tau \circ \Phi_t^{-1}(\etab)] e^{k_y(\tau-t)} d\tau,
\end{equation}
\begin{equation} \label{DDh:explicit}
\begin{array}{ll}
\DDh[t,\etab]=
&
\Dh^0[\Phi_t^{-1}(\etab)] e^{-k_D t}
\\ &
+ \int_0^t \left(
k_D D[\tau,\Phi_\tau \circ \Phi_t^{-1}(\etab)] -(\vm-\Pvh)(\tau).[\Phi_\tau \circ \Phi_t^{-1}](\etab)
\right) e^{k_D(\tau-t)}  d\tau.
\end{array}
\end{equation}
Thus, $\pwh$ and $\pvh$ are also explicit:
\begin{equation} \label{Pwh:explicit}
\pwh(t)=\pwh^0 - \int_0^t \int_{\SSS^2}\left( \lambda_y(\Yh-y)(\nabla \Yh\times\etab) +\lambda_D(\DDh-D)(\nabla \DDh\times\etab) \right)d\sigma_\eta d\tau,
\end{equation}
\begin{equation}
\begin{array}{ll} \label{Pvh:explicit}
\pvh(t)= \pvh^0 &- \int_0^t \int_{\SSS^2}\left( \lambda_D(\DDh-D)\etab+\lambda_y(\Yh-y)(\frac{1}{D}\etab\times(\etab\times\nabla \Yh)\right.\\
&\qquad\quad\left.+\lambda_D(\DDh-D)(\frac{1}{D}\etab\times(\etab\times\nabla \DDh)  \right)d\sigma_\eta d\tau.
\end{array}
\end{equation}

Now, let us introduce the map
$$\begin{array}{|cccl}
F_T: & E_T         & \rightarrow & C^0([0,T],\mathbb{R}^6) \\
   & (\Pwh,\Pvh) & \mapsto     & (\pwh,\pvh).
\end{array}$$
\\

\textbf{Step 1: Let us prove the existence of
$$T^*=T^*(M, R, \|(\wm,\vm)\|_{L^\infty(0,1)}, D_*, \|(D,y)\|_{C^1([0,1]\times\SSS^2)},k_y,k_D,\lambda_y,\lambda_D) \in (0,1)$$
such that $F_T$ maps $E_T$ into $E_T$ for every $T<T^*$.}
From now on, we assume that $T \in (0,1)$. Let $(\Pwh,\Pvh) \in  E_T$.
Thanks to the explicit expressions (\ref{Yh:explicit}), (\ref{DDh:explicit}) and the bounds (\ref{bound:diff_flow}), (\ref{yDC1b})
there exists a constant
$$C_1=C_1(M, R,\|(\wm,\vm)\|_{L^\infty(0,1)},D_*,\|(y,D)\|_{C^1([0,1]\times\SSS^2)},k_y,k_D)>0$$
such that
\begin{equation}
|\Yh(t,\etab)|, |\nabla \Yh(t,\etab)|, |\DDh(t,\etab)|, |\nabla \DDh(t,\etab)| \leqslant C_1, \quad \forall (t,\etab) \in [0,T]\times\SSS^2.
\end{equation}
Then, we deduce from (\ref{Pwh:explicit}), (\ref{Pvh:explicit}) and (\ref{def:R}) that
$$|\pwh(t)-\pw|, |\pvh(t)-\pv| \leqslant R + C_2 T, \forall t \in [0,T],$$
where
$$C_2=C_2(M, R, \|(\wm,\vm)\|_{L^\infty(0,1)}, D_*, \|(D,y)\|_{C^1([0,1]\times\SSS^2)},k_y,k_D,\lambda_y,\lambda_D)>0.$$
Thus if $T<R/C_2$, then $(\pwh,\pvh)$ belongs to $E_T$.
This ends Step 1 with $T^*:=\min\{1,R/C_2\}$.
\\

\textbf{Step 2: Let us prove the existence of
$$T^{**}=T^{**}(M,R,\|(\wm,\vm)\|_{L^\infty(0,1)},D_*,\|(y,D)\|_{C^1([0,1]\times\SSS^2)},k_y,k_D,\lambda_y,\lambda_D) \in (0,1)$$
such that $F_T$ is $(1/2)$-contractant on $E_T$ for every $T<T^{**}$.}
Let $(\Pwh^1,\Pvh^1)$, $(\Pwh^2,\Pvh^2) \in E_T$,
$\Phi^1_t$, $\Phi^2_t$ be the associated flows on $\SSS^2$
and $(\Yh^1, \DDh^1, \pwh^1, \pvh^1)$, $(\Yh^2, \DDh^2, \pwh^2, \pvh^2)$ be the associated solutions of (\ref{syst-fixedpoint}).
The regularity of the flow of (\ref{EDO}) with respect to the function $(\Pwh,\Pvh)$ justifies the existence of a constant
$C_3=C_3(R,\|(\wm,\vm)\|_{L^\infty(0,T)},D_*,\|D\|_{C^1([0,T]\times\SSS^2)})>0$
such that, for every $t \in [0,T]$,
$$\| \Phi_t^1-\Phi_t^2 \|_{C^1(\SSS^2)}, \| (\Phi_t^1)^{-1}-(\Phi_t^2)^{-1} \|_{C^1(\SSS^2)}
\leqslant C_3  \|(\Pwh^1-\Pwh^2,\Pvh^1-\Pvh^2)\|_{L^\infty(0,T)}\,.$$
Using the regularity of $y$, $\yh^0$, $D$, $\Dh^0$ and bound (\ref{borne_C1_yD0}), we deduce the existence of a constant
$C_4=C_4(M,R,\|(\wm,\vm)\|_{L^\infty(0,T)},D_*,\|(y,D)\|_{C^1([0,T]\times\SSS^2)},k_y,k_D)>0$
such that
$$\|(\Yh^1-\Yh^2,\DDh^1-\DDh^2)(t)\|_{C^1(\SSS^2)} \leqslant C_4 \|(\Pwh^1-\Pwh^2,\Pvh^1-\Pvh^2)\|_{L^\infty(0,T)}\,, \forall t \in [0,T]\,.$$
Thus, we deduce from (\ref{Pwh:explicit}) and (\ref{Pvh:explicit}) that
$$\| (\pwh^1-\pwh^2,\pvh^1-\pvh^2) \|_{L^\infty(0,T)} \leqslant C_5 T \|(\Pwh^1-\Pwh^2,\Pvh^1-\Pvh^2)\|_{L^\infty(0,T)}$$
for some constant $C_5$, that depends on the same quantitites as $C_4$.
This gives the conclusion with $T^{**}=1/(2C_5)$. $\hfill \Box$

\subsection{Bounds on solutions}

The goal of this section is the proof of the following result.

\begin{proposition} \label{Prop:borneC1_pr_iteration}
Let $R>0$. There exists
$$k_*=k_*(R,\|(\wm,\vm)\|_{L^\infty(0,1)}, D_*, \|(y,D)\|_{C^1([0,1]\times\SSS^2)})>0$$
such that for every $\pwh^0, \pvh^0 \in \mathbb{R}^3$ with (\ref{def:R}), $k_y, k_D>k_*$, $T^*>0$,
if $(\yh,\Dh,\pwh,\pvh) \in C^1([0,T^*] \times \SSS^2)^2 \times C^1([0,T^*],\mathbb{R}^3)^2$
is a solution of (\ref{Observer_CYpb}) on the time intervall $[0,T^*]$, then
$$\|(\yh,\Dh)(t)\|_{C^1(\SSS^2)} \leqslant \|(y,D)\|_{C^1(\SSS^2)} +1\,, \quad \forall t \in [0,T^*].$$
\end{proposition}

The proof of this proposition relies on the following technical Lemma.
\\

\begin{lemma}\label{lem:gradcont}
Let $M_0>0$. There exists $k_*=k_*(M_0)>0$ such that
\begin{itemize}
\item for every $\bold{a}\in C^1([0,+\infty)\times\SSS^2,\R^3)$, $b\in C^1([0,+\infty)\times\SSS^2,\R)$ such that
\begin{equation}\label{eq:a_tangent}
\bold{a}(t,\etab)\cdot\etab=0\,, \quad  \forall(t,\etab)\in[0,+\infty)\times\SSS^2\,,
\end{equation}
\begin{equation} \label{borneM0:a,b}
\|a(t,\etab)\|, \|b(t,\etab)\|, \left\|\dv{\bold{a}}{\etab}(t,\etab)\right\|,\left\|\dv{b}{\etab}(t,\etab)\right\|\leq M_0\,, \quad \forall(t,\etab)\in[0,+\infty)\times\SSS^2\,,
\end{equation}
\item for every $k>k_*$, $h^0 \in C^1(\SSS^2,\mathbb{R})$, $T^*>0$ and $h \in C^1([0,T^*] \times \SSS^2,\mathbb{R})$ solution on $[0,T^*]$ of
\begin{equation*}
 \left\{
\begin{array}{l}
\partial_t h (t,\etab)= \dv{h}{\etab}(t,\etab)\cdot \bold{a}(t,\etab)+ b(t,\etab) -kh(t,\etab)\,,\\
h(0,\etab)=h^0(\etab)\,,
\end{array}
\right.
\end{equation*}
\end{itemize}
then,
$$\|h(t,.)\|_{C^1(\SSS^2)} \leqslant \|h^0\|_{C^1(\SSS^2)} + 1\,, \forall t \in [0,T^*].$$
\end{lemma}

{\em Proof of Lemma \ref{lem:gradcont}:}
Let us consider the flow $\varphi(t,\etab)$ associated to the following equation
\begin{equation}
 \left\{
    \begin{array}{l}
\partial_t \varphi=-\bold{a}(t,\varphi)\\
\varphi(0,\etab)=\etab.
    \end{array}
\right.
\label{eq:varphi}
\end{equation}
The ODE has a local (in time) solution for any $\etab\in\SSS^2$ thanks to Cauchy-Lipschitz theorem because $\bold{a}$ is continuous in $t$ and $C^1$ in $\etab$. This solution lives on $\SSS^2$ thanks to assumption \eqref{eq:a_tangent}, thus no explosion is possible and $\varphi(t,\etab)$ is defined for every $(t,\etab)\in[0,+\infty)\times\SSS^2$.

The new function:
\begin{equation} \label{def:z_lemme1}
z(t,\etab):=h(t,\varphi(t,\etab))
\end{equation}
solves $\partial_t z(t,\etab)= b(t,\varphi(t,\etab))-kz$.
Thus
\begin{equation}\label{eq:lemme1}
z(t,\etab)=e^{-kt}h_0(\etab)+\int_0^t b(\tau,\varphi(\tau,\etab))e^{k(\tau-t)}d\tau\,.
\end{equation}

\textbf{Step 1: Bound on $h$.} For every $t \in [0,T^*]$, $\varphi(t,.)$ is a bijection of $\SSS^2$ thus
$$\begin{array}{ll}
\|h(t,.)\|_{L^\infty(\SSS^2)} & = \|z(t,.)\|_{L^\infty(\SSS^2)} \\
                              & \leqslant \|h^0\|_{L^\infty(\SSS^2)} e^{-kt} + \int_0^t M_0 e^{-k(t-\tau)} d\tau \quad \text{ by } (\ref{eq:lemme1}) \text{ and } (\ref{borneM0:a,b}) \\
                              & \leqslant \|h^0\|_{L^\infty(\SSS^2)} + \frac{M_0}{k} \\
                              & \leqslant \|h^0\|_{L^\infty(\SSS^2)} + 1 \quad \text{ when } k>M_0\,.
\end{array}$$

\textbf{Step 2: Bound on $\partial h/\partial \etab$.} From (\ref{def:z_lemme1}) and (\ref{eq:lemme1}), we deduce that
\begin{equation}\label{eq:lemme}
\begin{array}{lll}
\dv{h}{\etab}(t,\varphi(t,\etab))
& = & \dv{z}{\etab}(t,\etab) \left(\dv{\varphi}{\etab}(t,\etab)\right)^{-1} \\
& = & e^{-kt} \dv{h_0}{\etab}(\etab) \left(\dv{\varphi}{\etab}(t,\etab)\right)^{-1}\\
&   & + \int_0^t e^{-k(t-\tau)} \dv{b}{\etab}(\tau,\varphi(\tau,\etab)) \bold{\Sigma}(\tau,t,\etab) d\tau
\end{array}
\end{equation}
where
$$\bold{\Sigma}(\tau,t,\etab):=\dv{\varphi}{\etab}(\tau,\etab)\cdot\left(\dv{\varphi}{\etab}(t,\etab)\right)^{-1}\,.$$
Note that
$$\left\lbrace\begin{array}{l}
\dv{}{t}\left(\dv{\varphi}{\etab}(t,\etab)\right)^{-1}=\left(\dv{\varphi}{\etab}(t,\etab)\right)^{-1} \dv{\bold{a}}{\etab}(t,\varphi(t,\etab))\,,\\
\left(\dv{\varphi}{\etab}(0,\etab)\right)^{-1}=\bold{Id}
\end{array}\right.$$
and
$$\left\lbrace\begin{array}{l}
\dv{\bold{\Sigma}}{\tau}(\tau,t,\etab)=-\dv{\bold{a}}{\etab}(\tau,\varphi(\tau,\etab))\bold{\Sigma}(\tau,t,\etab)\,,\\
\bold{\Sigma}(t,t,\etab)=\bold{Id}\,.
\end{array}\right.$$
Thus, by Gronwall Lemma,
\begin{equation} \label{borne:phi_sigma}
\begin{array}{c}
\displaystyle
\left\| \left(\dv{\varphi}{\etab}(t,\etab)\right)^{-1} \right\| \leqslant e^{M_0 t}\,, \quad \forall (t,\etab) \in [0,T^*] \times \SSS^2\,,\\
\displaystyle
\| \bold{\Sigma}(\tau,t,\etab) \| \leqslant e^{M_0|t-\tau|}\,, \quad \forall (\tau,t,\etab) \in [0,T^*]^2 \times \SSS^2\,.
\end{array}
\end{equation}
We deduce from (\ref{eq:lemme}) and (\ref{borne:phi_sigma}) that
$$\begin{array}{ll}
\|\dv{h}{\etab}(t,\varphi(t,\etab))\|
& \leqslant  e^{(M_0-k)t} \|\dv{h_0}{\etab}\|_{L^\infty(\SSS^2)} + \int_0^t M_0 e^{(M_0-k)(t-\tau)} d\tau \\
& \leqslant \|\dv{h_0}{\etab}\|_{L^\infty(\SSS^2)} + \frac{M_0}{k-M_0} \quad \text{ when } k>M_0 \\
& \leqslant  \|\dv{h_0}{\etab}\|_{L^\infty(\SSS^2)} + 1 \quad \text{ when } k>2M_0\,. \\
\end{array}$$

This ends the proof of Lemma \ref{lem:gradcont} with $k_*=2M_0$. $\hfill \Box$
\\

{\em Proof of Proposition \ref{Prop:borneC1_pr_iteration}:} We apply Lemma \ref{lem:gradcont} to
$\widetilde{y}=\yh-y$ and $\widetilde{D}=\Dh-D$ with (see equations (\ref{eq:ytilde}) and (\ref{eq:gamtilde}))
$$a(t,\etab):=-\etab \times \left( (\wm-\pwh)(t) +\frac{1}{D(t,\etab)}\etab\times (\vm-\pvh)(t) \right) $$
$$b_1(t,\etab):= \nabla y(t,\etab)\cdot\left(\etab \times(\pwh(t)-\pw)  +\frac{1}{D(t,\etab)}\etab\times (\pvh(t)-\pv)  \right) $$
$$b_2(t,\etab):=  \nabla D(t,\etab)\cdot\left(\etab \times(\pwh(t)-\pw)  +\frac{1}{D(t,\etab)}\etab\times (\pvh(t)-\pv)  \right) + \pvh(t)\cdot\etab - \pv\cdot\etab\,.$$
Note that $\|(\pvh,\pwh)(t)\|$ is bounded uniformly with respect to $t \in [0,T^*]$ by a constant
that depends only on $R$ thanks to the decrease of the Lyapunov function $V$.
Thus, the assumptions of  Lemma \ref{lem:gradcont} are satisfied with
$$M_0=M_0(R,\|(\wm,\vm)\|_{L^\infty(0,1)}, D_*, \|(y,D)\|_{C^1([0,1]\times\SSS^2)})>0\,.$$
Therefore, there exists a constant $k_*>0$ (that depends on the same quantities as $M_0$)
such that, for every $k_y, k_D>k_*$ then
$$\|(\yh-y,\Dh-D)(t)\|_{C^1(\SSS^2)} \leqslant 1\,, \quad \forall t \in [0,T^*]\,,$$
which gives the conclusion. \hfill $\Box$

\subsection{Global solutions}

The goal of this section is the proof of Proposition \ref{Prop:WP1} thanks to Propositions \ref{Prop:WP1_local} and \ref{Prop:borneC1_pr_iteration}.

Let $\pwh^0, \pvh^0 \in \mathbb{R}^3$ and $R := \| (\pwh^0-\pw, \pvh^0-\pv)\|$.
We define
$$M:=\|(y,D)\|_{C^1(\SSS^2)} +1\,.$$
By Proposition \ref{Prop:WP1_local}, there exists a time $T=T(M,R)>0$ and a unique local solution of (\ref{Observer_CYpb}) defined on $[0,T]$.
By Proposition \ref{Prop:borneC1_pr_iteration}, we have
$$\|(\yh,\Dh)(T)\|_{C^1(\SSS^2)} \leqslant M\,.$$
By decreasing of the Lyapunov function, we have
$$\| (\pwh-\pw,\pvh-\pv)(T) \| \leqslant \sqrt{2 V(T)} \leqslant \sqrt{2 V(0)} = \| (\pwh^0-\pw,\pvh^0-\pv) \| = R$$
Thus we can apply Proposition \ref{Prop:WP1_local} with initial condition at $t=T$ and
we get a solution defined on $[0,2T]$. By iterating this argument, we get a solution defined for every $t \in [0,+\infty)$. \hfill $\Box$

\section{Appendix: Continuity of the flow}
\label{sec:app_cont}

The goal of this section is to prove Proposition \ref{Prop:WP3}.
We use the same notations as in the proof of Proposition \ref{Prop:WP1_local}.

Let $R \geqslant \| (\pwh^n(0)-\pw, \pvh^n(0)-\pv)\|$ for every $n \in \mathbb{N}$.
Thanks to the convergences in assumption, there exists $T_1>0$ such that the maps
$(F^n_{T_1})_{n \in \mathbb{N}}$ and $F_{T_1}$ associated to $(\vit^n,\om^n,y^n,D^n)_{n \in \mathbb{N}}$ and $(\vit,\om,y,D)$
are $(1/2)$-contractions of the same space $E_{T_1}$.

\textbf{Step 1: We prove that, for every $(\Pwh,\Pvh) \in E_{T_1}$ then
$$\|(F^n_{T_1}-F_{T_1})(\Pwh,\Pvh)\|_{L^\infty(0,T_1)}  \rightarrow 0 \text{ when } n \rightarrow + \infty.$$} Let $(\Pwh,\Pvh) \in E_{T_1}$. The first component of $(F^n_{T_1}-F_{T_1})(\Pwh,\Pvh)$ is
$$(F^n_{T_1}-F_{T_1})(\Pwh,\Pvh)^{(1)}(t)=(\pwh^n-\pwh)(0) + A_n(t) + B_n(t)$$ where, for every $t \in [0,T_1]$ (integrations by part)
$$A_n(t)=-\int_0^t \int_{\SSS^2} \lambda_y  \Big(
 (\Yh^n-y^n-\Yh+y)(\nabla \Yh^n \times \etab)   - (\Yh^n-\Yh) \nabla(\Yh-y) \times \etab \Big) d\sigma_\eta d\tau,$$
$$B_n(t)=-\int_0^t \int_{\SSS^2} \lambda_D  \Big(
(\DDh^n-D^n-\DDh+D)(\nabla \DDh^n \times \etab) - (\DDh^n-\DDh)  \nabla (\DDh-D)  \times \etab
\Big) d\sigma_\eta d\tau.$$
For every $t \in [0,T_1]$, we have
\begin{equation} \label{decA}
\begin{array}{ll}
|A_n(t)| \leqslant \lambda_y T_1 C \Big(
&
 \|\Yh^n-y^n-\Yh+y\|_{L^\infty((0,T_1)\times\SSS^2)} \|\nabla \Yh^n\|_{L^\infty((0,T_1)\times\SSS^2)}
\\ &
+ \|\nabla (\Yh-y)\|_{L^\infty((0,T_1)\times\SSS^2)} \|\Yh^n-\Yh \|_{L^\infty((0,T_1)\times\SSS^2)}
\Big).
\end{array}
\end{equation}
Thanks to (\ref{hyp:CV(v,w)}), we also have $\Phi_t^n \rightarrow \Phi_t$ in $C^0([0,T_1]\times\SSS^2,\SSS^2)$ (consequence of Gronwall lemma).
From the explicit expressions of $\Yh^n, \Yh, \DDh^n, \DDh$ (see \eqref{Yh:explicit} and \eqref{DDh:explicit}) and the convergence (\ref{hyp:CV(y,G)}),
we deduce that
\begin{equation} \label{CV:(Yh,GGh)}
(\Yh^n,\DDh^n) \rightarrow (\Yh,\DDh) \quad \text{ in } C^0([0,T_1]\times \SSS^2,\mathbb{R})^2.
\end{equation}
Moreover, thanks to (\ref{hyp:bound(y,G)}), the quantities $\|\nabla \DDh^n\|_{L^\infty((0,T_1)\times\SSS^2)}$,$\|\nabla \Yh^n\|_{L^\infty((0,T_1)\times\SSS^2)}$ are uniformly bounded.
Thus, the right hand side of (\ref{decA}) converges to zero as $n \rightarrow + \infty$.
We may work similarly to prove that $B_n(t)$ and $(F^n_{T_1}-F_{T_1})(\Pwh,\Pvh)^{(2)}(t)$ converge to zero uniformly with respect to $t \in [0,T_1]$.
This ends the proof of the first step.

\textbf{Step 2: Conclusion.}
We have
\begin{equation}
\begin{array}{ll}
          & \| (\pwh^n-\pwh,\pvh^n-\pvh) \|_{L^\infty(0,T_1)} \\
=         & \| F^n_{T_1}(\pwh^n,\pvh^n) - F_{T_1}(\pwh,\pvh)  \|_{L^\infty(0,T_1)} \\
\leqslant & \|F^n_{T_1}(\pwh^n,\pvh^n)-F^n_{T_1}(\pwh,\pvh)\|_{L^\infty(0,T_1)} +\|(F^n_{T_1}-F_{T_1})(\pwh,\pvh)\|_{L^\infty(0,T_1)}
\end{array}
\end{equation}
and $F^n_{T_1}$ is a $(1/2)$-contraction of $E_{T_1}$ thus
$$\| (\pwh^n-\pwh,\pvh^n-\pvh) \|_{L^\infty(0,T_1)} \leqslant 2 \|(F^n_{T_1}-F_{T_1})(\pwh,\pvh)\|_{L^\infty(0,T_1)} \rightarrow 0$$
thanks to the first step. Iterating this argument on a finite number of intervals $[T_1,2T_1]$, $[2T_1,3T_1]$, etc, we conclude that (\ref{CV(pwh,pvh)}) holds.
Let $\varphi_t$ be the flow associated to the equation
$$\frac{d\etab}{dt}=\etab \times \Big( \wm(t)-\pwh(t) +\frac{1}{D(t,\etab)}\etab\times [\vm(t)-\pvh(t)] \Big)$$
and $\varphi_t^n$ similarly associated to $(\pwh^n,\pvh^n)$.
Then, by (\ref{CV(pwh,pvh)}), $\varphi_t^n \rightarrow \varphi_t$ in $C^0([0,T]\times \SSS^2,\SSS^2)$.
Thanks to the explicit expressions of $(\yh^n\Dh^n)$ and $(\yh,\Dh)$ in terms of $\varphi_t^n$ and $\varphi_t$
(obtained by the characteristic method, as (\ref{Yh:explicit}) and (\ref{DDh:explicit})),
we deduce that (\ref{CV(yh,Gh)}) holds. This ends the proof of Proposition \ref{Prop:WP3}. \hfill $\Box$

\section{Appendix: Basic formulae of differential geometry on $\SSS^2$: $\bold{P}$ is a constant vector}
\subsection{$\nabla(\etab\cdot \bold{P})=-\etab\times(\etab\times \bold{P})$: proof of formula \eqref{formule:gradient_scalaire}}
\label{subsubsec:gradient}
By definition of the differential of a mapping defined on $\SSS^2$,
\begin{equation}
\nabla(\etab\cdot \bold{P})\cdot \delta\etab=\delta(\etab\cdot \bold{P}) \quad\forall \delta\etab\quad s.t. \quad\delta\etab\cdot\etab=0
\end{equation}
or
\begin{equation}
\nabla(\etab\cdot \bold{P})\cdot \delta\etab=\delta\etab\cdot \bold{P}
\label{eq:difdugradient}
\end{equation}
On the other hand, $\nabla(\etab\cdot \bold{P})\cdot\etab=0$ because the gradient in $\etab$ lives in the tangent plane of $\SSS^2$ in $\etab$.
Since $\delta \etab$ is in the tangent plane, it yields that $\nabla(\etab\cdot \bold{P})$ can be identified as the projection of $\bold{P}$ in the tangent plane. This writes $\bold{P}-(\etab\cdot \bold{P})\etab$, or equivalently $-\etab\times(\etab\times \bold{P})$.

\subsection{$\Delta(\etab\cdot \bold{P})=-2\eta\cdot \bold{P}$: proof of formula \eqref{formule:laplacien}}
\label{subsubsec:laplacien}
By definition, the Laplacian of a function is the divergence of its gradient. Thus $\Delta(\etab\cdot \bold{P})=\nabla\cdot (-\etab\times(\etab\times \bold{P}))$. Let us define $\bold{F}(\etab)=\etab\times(\etab\times \bold{P})$. The divergence of a vector field is the derivative of the area of a surface element as it evolves along the flow defined by the vector field. More precisely, as $\delta \etab_1$ and $\delta \etab_2$ are two infinitesimal vectors of $T_\eta$ for $\etab \in \SSS^2$, they define a surface element whose area is  $(\delta \etab_1 \times \delta \etab_2)\cdot \etab$. As $\etab$, $\delta \etab_1$ and $\delta \etab_2$ evolve according to
\begin{equation}
\dv{\etab}{s}=F(\etab)
\label{eq:flotEta}
\end{equation}
the divergence verifies:
\begin{equation}
\dv{(\delta \etab_1 \times \delta \etab_2)\cdot \etab}{s}=\nabla\cdot F(\etab)(\delta \etab_1 \times \delta \etab_2)\cdot \etab
\label{eq:defDiv}
\end{equation}
It is developped as
\begin{align*}
\dv{(\delta \etab_1 \times \delta \etab_2)\cdot \etab}{s}=&(\dv{\delta \etab_1}{s} \times \delta \etab_2)\cdot \etab+(\delta \etab_1 \times \dv{\delta \etab_2}{s})\cdot \etab+(\delta \etab_1 \times \delta \etab_2)\cdot \dv{\etab}{s}
\end{align*}
Using \eqref{eq:flotEta} and the definition of $F(\etab)$ applied to $\delta \etab_1$, $\delta \etab_2$ and $\etab$:
\begin{align*}
\dv{(\delta \etab_1 \times \delta \etab_2)\cdot \etab}{s}=&((-\delta\etab_1\times(\etab\times \bold{P})-\etab\times(\delta\etab_1\times \bold{P})) \times \delta \etab_2)\cdot \etab\\
&+(\delta \etab_1 \times (-\delta\etab_2\times(\etab\times \bold{P})-\etab\times(\delta\etab_2\times \bold{P})))\cdot \etab\\
&+(\delta \etab_1 \times \delta \etab_2)\cdot (-\etab\times(\etab\times \bold{P}))
\end{align*}
Rearranging the terms, and since  $\delta\etab_1$ and $\delta\etab_2 \in T_\eta$ yields:
\begin{align*}
\dv{(\delta \etab_1 \times \delta \etab_2)\cdot \etab}{s}=-2(\etab\cdot \bold{P})(\delta \etab_1 \times \delta\etab_2)\cdot \etab
\end{align*}
This concludes the proof of the formula by the definition \eqref{eq:defDiv} of the divergence.

\subsection{$\nabla\cdot(\etab\times \bold{P})=0$: proof of formula \eqref{formule:divergence}}
\label{subsubsec:rotation}
It is obvious as soon as one realize that $(\etab\times \bold{P})$ is a rotation, thus the lengths are left unchanged by this transformation.



\end{document}